\documentclass{article}
\usepackage[greek,american]{babel}
\usepackage{amssymb}
\usepackage{amsfonts}
\usepackage[dvips]{graphicx}
\begin{document}
\newtheorem{remark}{Remark}[section]
\newtheorem{pigsfly1}{Definition}[section]
\newtheorem{rigor1}{Theorem}[section]
\newtheorem{rigor2}{Proposition}[section]
\newtheorem{smallTheorem}{Lemma}[section]
\newtheorem{help1}{Example}[section]
\title{Asymptotic Behavior of Solutions of Complex Discrete Evolution Equations: The Discrete Ginzburg-Landau Equation}

\author{Nikos I. Karachalios$^{\dag}$, Hector E. Nistazakis$^{*}$, Athanasios N. Yannacopoulos$^{\ddag}$}
\maketitle
\begin{abstract}
We study the asymptotic behavior of complex discrete evolution equations of Ginzburg-
Landau type. Depending on the nonlinearity and the data of the problem, we find different 
dynamical behavior ranging from global existence of solutions and global attractors, to 
blow up in finite time. We provide estimates for the blow up time, depending not only on 
the initial data but also on the size of the lattice. Some of the theoretical results, are tested  
by numerical simulations.  
\end{abstract}
\section{Introduction}
Several theoretical and experimental studies  performed on spatially discrete systems, have proved that such systems display very reach dynamical behavior, even in the one-dimensional space. The Discrete Complex Ginzburg-Landau Equation (DCGL),
\begin{eqnarray}
\label{introlat1}
\dot{u}_n-(\lambda+i\alpha)(u_{n-1}-2u_n+u_{n+1})=(k+i\beta)|u_n|^2u_n+\gamma u_n,
\end{eqnarray}
(where the $n$-index ranges over the $1D$-lattice), is a particular discrete evolution equation, whose dynamics can lead to extraodinary complicated behavior, ranging from spatiotemporal intermittency and dispersive chaos, to self-localization phenomena and the existence of discrete solitons. One of the most interesting applications where the DCGL equation (\ref{introlat1}) may arise, is in the description of the evolution of Taylor and frustrated vortices, in hydrodynamic systems of low dimensionality, and it has been proved a fairly reasonable mathematical  model for investigating phenomena related to weak turbulence \cite{PLeGal,WilCa1}. The famous Discete Nonlinear Schr\"odinger Equation (DNLS) (obtained from (\ref{introlat1}), in the case $\lambda=k=0$) \cite{Camp, EilJo, Kevrekidis},  is encountered in several diverse branches of physics, ranging  from supeconductivity and nonlinear optics, to the Bose-Einstein condensates.

The aim of this work, is to provide some results, concerning the asymptotic behavior of solutions of discrete evolution equations of the form, 
\begin{eqnarray}
\label{lat1}
\dot{u}_n-(\lambda+i\alpha)(u_{n-1}&-&2u_n+u_{n+1})=F(u_n)+\gamma u_n,\\
\label{lat2}
u_n(0)&=&u_{n,0},
\end{eqnarray}
where $\lambda\geq 0$ and $\alpha,\gamma\in\mathbb{R}$, and the classical cubic nonlinearity of (\ref{introlat1}), has been replaced by more general nonlinear interactions. The lattice may be infinite $(n\in\mathbb{Z})$ or finite $(|n|\leq N)$, supplemented with Dirichlet boundary conditions. 
Mainly we are interested in nonlinear interactions of the form $F(s)=(k+i\beta)|s|^p$ (non-gauge interaction), or $F(s)=(k+i\beta)|s|^{p-1}s$, for some $p>1$ (gauge interaction). The case $\alpha=\beta=0$, corresponds to the discrete Ginzburg-Landau equation with real coefficients (DRGL).

Depending on the type of the nonlinearity, the length of the lattice and  the``size'' of the initial data, the dynamical behavior of solutions ranges from  the existence of finite time singularities (blow-up in finite time) to the existence of global attractors in appropriate and physically justified phase spaces. 
 
With respect to  the issue of global existence and blow-up of solutions, it has been observed numerically that discreteness may have important effects. For example, in the case of the conservative DNLS with gauge nonlinear interaction, solutions exist globally, independently of the choice of the initial data and the strength of the nonlinearity. This is in contrast with the NLS continuous counterpart, for which solutions may blow-up in finite time. 
As it is shown in Section 3, the DCGL and DRGL equations, serve as a discrete models, whose behavior differs with respect to this issue, since both in the case of the gauge and non-gauge interactions, solutions may blow-up in finite time, at least in the case of a finite lattice. On the other hand, this behavior is in agreement with that of the CGL and RGL partial differential equation (see the recent results of \cite{Ozawa}). 

In the case of a dissipative gauge nonlinearity the dynamics of the finite system can be described by a global attractor. A simple but interesting observation, is that the entry time to the absorbing ball is independent of the size of the initial data, another diferrence with the CGL partial differential equation. To test the theoretical estimates for the blow-up time, we perform a first attempt on the derivation of some numerical simulations, for the case of non-gauge nonlinearity.

Motivated by the pioneering work of \cite{Bates} for lattice dynamical systems of first order and \cite{SZ2} for extensions to various multidimensional lattices of first and second order, we devote Section 4 to the study of (\ref{lat1})-(\ref{lat2}) in the case of an infinite lattice, as an infinite dimensional dynamical system.  However, our approach differs from the aforementioned references, since we study the DCGL in {\em  a particular class of weighted sequence spaces}. These spaces, cover (but are not limited to) the case of {\em exponential localized} solutions for lattice differential equations. Substantial physical motivation is given for the study of (\ref{lat1})-(\ref{lat2}) in such spaces, since many of the physically interesting solutions of complex discrete evolution equations, present strong spatial localization properties. This is true for instance for soliton solutions or breathers \cite{Mackayexp,Camp, EilJo, FlachWillis, Kevrekidis, MackayAubry, Kevrekidis2}.

As in the case of the finite lattice, in the case of a dissipative nonlinearity we prove the existence of a global attractor in the weighted phase space. Let us note once again (and in connection with the corresponding result for the finite lattice), that discreteness enforces us to impose some restrictions on the parameters, which may differ from the usual restrictions on the dissipation parameter, appearing in the case of the continuous CGL equation. We conclude by discussing the approximation of the global attractor of spatial localized solutions, by  the
dynamics of the finite dimensional system. 
\section{Phase spaces and local existence of solutions}
This preliminary section is devoted to the definition of the appropriate phase spaces and to the basic results on local existence of solutions for (\ref{lat1})-(\ref{lat2}).

The case of the finite lattice for (\ref{lat1})-(\ref{lat2}) $(|n|\leq N)$, with Dirichlet boundary conditions,  will be considered in the finite dimensional Hilbert space
$\mathbb{C}^{2N+1}$ endowed with the usual inner product and Euclidean norm,
\begin{eqnarray}
\label{discn}
(\phi,\psi)_2:=\mathrm{Re}\sum_{n=-N}^{n=N}\phi_n\overline{\psi_n},\;\;||\psi||^2_2:=\sum_{n=-N}^{n=N}|\psi_n|^2,\;\;\phi,\,\psi\in \mathbb{C}^{2N+1},
\end{eqnarray}
We consider now the operators
\begin{eqnarray}
\label{finself}
(A_d\psi)_{|n|\leq N}:=\psi_{n-1}-2\psi_n+\psi_{n+1},\;\;(B_d\psi)_{|n|\leq N}=\psi_{n+1}-\psi_n.
\end{eqnarray}
It can be easily checked (see also \cite[pg. 117]{Akriv}) that
\begin{eqnarray}
\label{byparts}
(-A_d\psi,\psi)_2=\sum_{n=-N}^{n=N}|\psi_{n+1}-\psi_n|^2,\;\;(-A_d\phi,\psi)_2=(B_d\phi,B_d\psi)_2.
\end{eqnarray}
Hence, we may also consider the inner product and the corresponding norm in $\mathbb{C}^{2N+1}$, 
\begin{eqnarray}
\label{discSob}
(\phi,\psi)_{1,2}:=(B_d\phi,B_d\psi)_2+(\phi,\psi)_2,\;\;||\psi||_{1,2}:=\sum_{n=-N}^{n=N}(|\psi_{n+1}-\psi_n|^2+|\psi_n|^2).
\end{eqnarray}
In our analysis, we shall also use for any $1\leq p\leq\infty$ the norms
\begin{eqnarray*}
||\psi||_{p}=\left(\sum_{n=-N}^{n=N}|\psi_n|^p\right)^{1/p},\;\;||\psi||_{\infty}=\max\left\{|\psi_n|\;:\;|n|\leq N\right\},\;\;\psi\in \mathbb{C}^{2N+1}.
\end{eqnarray*}
For any $1\leq p\leq q\leq\infty$, there exist constants $c_1,c_2$ depending on $N$, 
\begin{eqnarray}
\label{otherfnorms}
c_1||\psi||_{p}\leq ||\psi||_{q}\leq c_2||\psi||_{p},\;\;\psi\in \mathbb{C}^{2N+1}.
\end{eqnarray}
We note that the norm in (\ref{discSob}),  is equivalent with the norm in (\ref{discn}). 

For the case of an infinite lattice $(n\in\mathbb{Z})$, a first natural choice for the phase space is to consider complexifications of the usual real sequence spaces, denoted by
\begin{eqnarray}
\label{lp}
{\ell}^p:=\left\{u=(u_n)_{n\in\mathbb{Z}}\in\mathbb{C}\;\;:\;\;
||u||_{\ell^p}:=\left(\sum_{n\in\mathbb{Z}}|u_n|^p\right)^{\frac{1}{p}}<\infty\right\}.
\end{eqnarray}
Between $\ell^p$, spaces the following elementary embedding relation  \cite[pg. 145]{HiLa} holds,
\begin{eqnarray}
\label{lp1}
\ell^q\subset\ell^p,\;\;\;\; ||u||_{\ell^p}\leq ||u||_{\ell^q}\,\;\; 1\leq q\leq p\leq\infty.
\end{eqnarray}
For $p=2$ we get the usual Hilbert space of square-summable (complex) sequences endowed with the real scalar product
\begin{eqnarray}
\label{lp2}
(u,v)_{\ell^2}=\mathrm{Re}\sum_{{n\in\mathbb{Z}}}u_n\overline{v_n},\;\;u,\,v\in\ell^2.
\end{eqnarray}
Of particular interest is also the existence result in $\ell^1$ which  can be considered as the space of ``discrete regularity'', in the sense suggested by (\ref{lp1}).

To cover the  situation of spatially localized solutions, we study (\ref{lat1})-(\ref{lat2}) in weighted spaces, with properly chosen weight functions. We consider a 
weight function $\theta_n$, which is an increasing function of $|n|$, satisfying for all $n\in\mathbb{Z}$, the following condition: there exist constants $D,\underline{d},\overline{d}>0$, such that
\begin{eqnarray}
\mathrm{(WS)} \left \{ \begin{array}{ccc}
&1\leq \theta_n\nonumber \\
&\mid \theta_{n+1}-\theta_{n} \mid \le D \theta_n
\nonumber \\
&\underline{d} \theta_n \le \theta_{n+1}\le \overline{d}\theta_{n},
\nonumber
\end{array}
\right .
\end{eqnarray}
and we introduce the weighted  spaces $\ell^p_\theta$, $$\ell_\theta^p=\{u_n \in {\mathbb C}:\;
\mid\mid u \mid\mid_{\ell_\theta^p}^p:=\sum_{n\in\mathbb{Z}} \theta_n \mid
u_n \mid^p < \infty\}.$$  It can easily be seen that the space
$\ell_\theta^2$ is a Hilbert space, endowed with the norm $\mid\mid \;
\mid\mid_{\ell_\theta^2}$ and scalar product
\begin{eqnarray}
\label{weightscal}
(u,v)_{\theta}=\mathrm{Re}\sum_{n\in\mathbb{Z}}\theta_nu_n\overline{v_n},\;\;u,\,v\in\ell^2_\theta.
\end{eqnarray}
Such spaces are the discrete analogue of
weighted $L^{p}$ spaces. It follows from $\mathrm{(WS)}$, that for any $1\leq p\leq\infty$,
\begin{eqnarray}
\label{we1}
\ell^p_\theta\subset\ell^p,\;\;\;||u||_{\ell^p}\leq ||u||_{\ell^p_{\theta}},\;\;\;1\leq p\leq\infty.
\end{eqnarray}
Moreover, we observe by using (\ref{lp1}) and $\mathrm{(WS)}$, that for any $1\leq q\leq p\leq\infty$, 
 similar embedding relations to (\ref{lp1}), hold for the weighted sequence spaces, that is
\begin{eqnarray}
\label{wlp1}
\ell^q_{\theta}\subset\ell^p_{\theta},\;\;\;\; ||u||_{\ell^p_{\theta}}\leq ||u||_{\ell^q_{\theta}}\,\;\; 1\leq q\leq p\leq\infty.
\end{eqnarray}
Let us remark that a choice for a weight function satisfying $\mathrm{(WS)}$, is the
exponential function $\theta_n=exp(\mu \mid n\mid)$ for
$\mu >0$. Existence of solutions in such spaces will provide
us with the existence of (exponentially) localized solutions for
the DGL equation. An instance where such spaces have been used is in
\cite{Mackayexp} where  the existence of exponentially localized
solutions has been studied in conservative lattices using a
continuation argument, related to the anti-integrable limit. 

For local existence os solutions, we shall examine the following examples of nonlinearities $F:\mathbb{C}\rightarrow\mathbb{C}$: \vspace{.2cm}
\newline
$\mathrm{(N_1)}$\ \ $F(0)=0$ and there exist constants $c>0$, $p>1$ such that
$|F(z_1)-F(z_2)|\leq c(|z_1|^{p-1}+|z_2|^{p-1})|z_1-z_2|$, or alternatively
\newline
$\mathrm{(N_2)}$\ \ $F(z)=f(|z|^2)z$ where $f,\,f':\mathbb{R}\rightarrow\mathbb{R}$, are continuous.
\vspace{.2cm}
\newline
First, we shall need some information on the nonlinear maps defined by the nonlinear interactions, provided by the following
\begin{smallTheorem}
\label{LipN}
Let $X$ be either the space $\mathbb{C}^{2N+1},\ell^2,\ell^2_\theta,\ell^1$ and assume that $F:\mathbb{C}\rightarrow\mathbb{C}$ satisfies $\mathrm{(N_1)}$ or $\mathrm{(N_2)}$. Then the function $F$ defines an operator (still denoted by $F$) $$F:X\rightarrow X,\;\;(F(u))_{n\in\mathbb{Z}}:=F(u_n),$$ which is Lipschitz continuous on bounded sets of $X$.
\end{smallTheorem}
{\bf Proof:}\ \ We  focus on the case of the infinite lattice, since the treatment of the finite lattice is almost the same. More precisely, we present only the case $X=\ell^2_{\theta}$ (since the case for $\ell^2$ and $\ell^1$ is similar). 
Let $u\in B_R$ a closed ball in $\ell^2_{\theta}$, of center $0$ and radius $R$.  We have from (\ref{wlp1}) that
\begin{eqnarray}
\label{ufp1}
||\mathrm{F}(u)||^2_{\ell^2_{\theta}}&\leq& c^2\sum_{n\in\mathbb{Z}}\theta_n|u_n|^{2p}
=c^2||u||_{\ell_{\theta}^{2p}}^{2p}\leq c^2||u||_{\ell^2_{\theta}}^{2p},
\end{eqnarray}
hence $\mathrm{F}:\ell^2_{\theta}\rightarrow\ell^2_{\theta}$,  is  bounded on bounded
sets of $\ell^2_{\theta}$.

For $u,v\in B_R$, we observe by using the embedding $\ell^2_{\theta}\subset\ell^{\infty}$ (provided by  (\ref{lp1}) and (\ref{we1})), that
\begin{eqnarray}
\label{ufp2}
||\mathrm{F}(u)-F(v)||^2_{\ell^2_{\theta}}&\leq& c^2\sum_{n\in\mathbb{Z}}\theta_n(|u_n|^{p-1}+|v_n|^{p-1})^2|u_n-v_n|^2\\
&\leq& c^2\sup_{n\in\mathbb{Z}}\left[(|u_n|^{p-1}+|v_n|^{p-1})^2\right]
\sum_{n\in\mathbb{Z}}\theta_n|u_n-v_n|^2\leq c^24R^{2(p-1)}||u-v||_{\ell^2_{\theta}}^2,
\end{eqnarray}
justifying that the map $F:\ell^2_{\theta}\rightarrow\ell^2_{\theta}$, is
Lipschitz continuous on bounded sets of $\ell^2_{\theta}$, with Lipschitz constant $L(R)=c2R^{(p-1)}$.\newline

For the case $\mathrm{(N_2)}$, we have
\begin{eqnarray}
\label{prop1}
||\mathrm{F}(u)||^2_{\ell^2_{\theta}}=\sum_{n\in\mathbb{Z}}\theta_n|f(|u_n|^2)|^2|u_n|^2.
\end{eqnarray}
Since $f:\mathbb{R}\rightarrow\mathbb{R}$ is continuous, there exists a monotone increasing $\mathrm{C^1}$-function
$g:\mathbb{R}\rightarrow\mathbb{R}$ such that
\begin{eqnarray}
\label{prop2}
|f(\rho)|\leq g(|\rho|),\;\;\mbox{for all}\;\;\rho\in\mathbb{R},
\end{eqnarray}
(see e.g \cite[p.g 796]{zei85}). Now by using (\ref{prop2}), we get
\begin{eqnarray}
\mid\mid F(u) \mid\mid_{\ell_{\theta}^{2}}^{2}&=&\sum_{n\in\mathbb{Z}}
\theta_n\mid f(\mid u_n \mid^2)\mid^2\mid u_n \mid^2 \le
\sum_{n\in\mathbb{Z}}\theta_n g(\mid u_n\mid^2)^2\mid u_n\mid^2
\nonumber\\
&\le& \sum_{n\in\mathbb{Z}}\theta_n g(\mid\mid u\mid\mid_{\ell_{\theta}^{2}}^2)^2\mid u_n\mid^2
\le\{\max_{\rho\in[0,R^2]} g(\rho)\}^2 \sum_{n\in\mathbb{Z}}\theta_n \mid u_n\mid^2 \le c(R)
\mid\mid u\mid\mid_{\ell_{\theta}^{2}}. \nonumber
\end{eqnarray}
for some positive constant $c(R)$. Thus we conclude, that the operator $F$ is bounded on bounded sets of $\ell_{\theta}^{2}$. To check the Lipschitz property, we may see, that for some $\nu \in (0,1)$, 
\begin{eqnarray}
\mid\mid F(u)-F(v)\mid\mid_{\ell_{\theta}^{2}} 
&\le& 2\sum_{n\in\mathbb{Z}} \theta_n  \mid f(\mid u_n\mid^2)\mid^2 \mid u_n -v_n\mid^{2}\nonumber\\
&&+ 2\sum_{n\in\mathbb{Z}}\theta_n \mid f^{'} (\nu \mid u_n\mid^2
+(1-\nu)\mid v_n\mid^2)\mid^{2}(\mid u_n\mid +\mid
v_n\mid)^{2}  \mid v_n\mid^{2} \mid u_n -v_n\mid^2.
\nonumber
\end{eqnarray}
For an appropriate $\mathrm{C^1}$-function
$g_1:\mathbb{R}\rightarrow\mathbb{R}$, we get the inequality
\begin{eqnarray}
\mid\mid F(u)-F(v)\mid\mid_{\ell_{\theta}^{2}}^{2} &\le& 2\sum_{n\in\mathbb{Z}} \theta_{n}\mid f(\mid\mid  u\mid\mid_{\ell_{\theta}^{2}}^{2})\mid^{2}\mid u_n -v_n\mid^{2}\nonumber\\
&&+2\{ \max_{\rho \in [0,2R^2]}
g_1(\rho) \}^{2} c(R) \sum_{n\in\mathbb{Z}}\theta_n \mid
u_{n}-v_n\mid^{2}
\nonumber \\
&\le& c(R)( \mid\mid u -v \mid\mid_{\ell_{\theta}^{2}}^{2}.
\nonumber
\end{eqnarray}

Now in the case of $\ell^1$, we consider as an example the case $\mathrm{(N_1)}$. For the Lipschitz condition we have
\begin{eqnarray*}
||\mathrm{F}(u)-F(v)||_{\ell^1}&\leq& c\sum_{n\in\mathbb{Z}}(|u_n|^{p-1}+|v_n|^{p-1})|u_n-v_n|\\
&\leq& c\sup_{n\in\mathbb{Z}}\left[(|u_n|^{p-1}+|v_n|^{p-1})\right]
\sum_{n\in\mathbb{Z}}|u_n-v_n|\leq c2R^{(p-1)}||u-v||_{\ell^1}.
\end{eqnarray*}
This concludes the proof of the Lemma.\ \ $\diamond$

For local existence of solutions, one could apply alternatively, a semigroup approach or existence theorems of ordinary differential equations on Banach spaces, depending on the choice of the phase space. The first approach could be of interest, as a starting point for further investigations on the properties of discrete operators \cite{Davies1}.
\subsection{Local existence in $\ell^2$}\ \ 
In what follows, a
complex Hilbert space $X$, endowed with the sesquilinear form
$B_{X}(\cdot ,\cdot)$ and the norm
$||\cdot||_{X}$, will be considered as a {\em real}
Hilbert space, endowed with the scalar product $(\cdot
,\cdot)_{X}=\mathrm{Re}\,B_{X}(\cdot
,\cdot)$. Let
$\mathbf{T}:D(\mathbf{T})\subseteq X\rightarrow X$, be  a $\mathbb{C}$-linear, self-adjoint, non-positive operator with dense domain $D(\mathbf{T})$, on the Hilbert space $X$, equipped with the scalar product $(\cdot ,\cdot)_{X}$. The space
$X_{\mathbf{T}}$, is the completion of $D(\mathbf{T})$ in the norm $||u||_{\mathbf{T}}^2=||u||^2_X-(\mathbf{T}u,u)_X$ for $u\in X_{\mathbf{T}}$. We denote
by $X_{\mathbf{T}}^*$ its dual, and by $\mathbf{T}^*$, the extension of $\mathbf{T}$ to the dual of $D(\mathbf{T})$, denoted by $D(\mathbf{T})^*$. 

For any $u,v\in\ell^2$
we consider  the linear operators $A,B,B^*:\ell^2\rightarrow\ell^2$,
\begin{eqnarray}
\label{diffop}
(Bu)_{n\in\mathbb{Z}}&=&u_{n+1}-u_{n},\;\;\;\;\;\;
(B^*u)_{n\in\mathbb{Z}}=u_{n-1}-u_{n},\\
(Au)_{n\in\mathbb{Z}}&=&(u_{n-1}-2u_n+u_{n+1}).
\end{eqnarray}
\begin{smallTheorem}
\label{semgen}
We assume that $\xi_1,\xi_2>0$. The operator $\mathbf{L}:\ell^2\rightarrow\ell^2$, $(\mathbf{L}u)_{n\in\mathbb{Z}}=(\xi_1+i\xi_2)(Au)_{n\in\mathbb{Z}}$ is the generator of a one parameter semigroup of $U(t)$ on $\ell^2$, that solves the underlying linear equation $\dot{u}=\mathbf{L}u$, namely $U(t)=\exp(\mathbf{L}t)$.
\end{smallTheorem}
{\it Proof:} We observe that the operator $\mathbf{L}$, is associated to a {\em non-symmetric} bilinear form on $\ell^2$ since for any $u=u_1+iu_2,v=v_1+iv_2\in\ell^2$,
\begin{eqnarray}
\label{sem1}
(\mathbf{L}u,v)_{\ell^2}=\xi_1\left\{(Bu_1,Bv_1)_{\ell^2}+(Bu_2,Bv_2)_{\ell^2}\right\}-\xi_2\left\{(Bu_2,Bv_1)_{\ell^2}-(Bu_1,Bv_2)_{\ell^2}\right\}.
\end{eqnarray}
However, we observe from (\ref{sem1}), that the operator $A_1u=\xi_1 Au$, satisfies the relations
\begin{eqnarray}
\label{lp6}
(A_1u,u)_{\ell^2}&=&-\xi_1||Bu||_{\ell^2}^2\leq 0,\\
\label{lp7}
(A_1u,v)_{\ell^2}&=&(u,A_1v)_{\ell^2},
\end{eqnarray}
therefore defines a self-adjoint operator on $D(A)=X=\ell^2$ and $A_1\leq 0$.  
We denote next by $\ell^2_{\xi_1}$, the Hilbert space with  the following scalar product and induced norm
\begin{eqnarray}
\label{lp3}
(u,v)_{\ell^2_{\xi_1}}&:=&\xi_1(Bu,Bv)_{\ell^2}+ (u,v)_{\ell^2},\\
\label{lp4}
||u||_{\ell^2_{\xi_1}}^2&:=&\xi_1||Bu||_{\ell^2}^2+||u||_{\ell^2}^2.
\end{eqnarray}
The usual norm of $\ell^2$ and (\ref{lp4}) are equivalent (see also \cite{SZ2}), since for some constant $c(\xi_1)>0$,
\begin{eqnarray}
\label{lp5}
||u||^2_{\ell^2}\leq ||u||^2_{\ell^2_{\xi_1}}\leq c||u||^2_{\ell^2}.
\end{eqnarray}
Note that the graph norm
$$||u||_{D(A_1)}=||A_1u||^2_{\ell^2}+||u||^2_{\ell^2},$$ is  an equivalent norm with the $\ell^2$-norm since
$$||u||_{\ell^2}^2\leq\xi_1^2\sum_{n\in\mathbb{Z}}|u_{n+1}-2u_n+u_{n-1}|^2+\sum_{n\in\mathbb{Z}}|u_n|^2\leq c||u||_{\ell^2}^2.$$
In our case, as it is indicated by
(\ref{lp3})-(\ref{lp6})-(\ref{lp7}), we may choose $X_{A_1}=\ell^2_{\xi_1}$
equipped with the norm $||u||_{A_1}^2=||u||_X^2-(A_1u,u)_X\equiv
||u||_{\ell^2_{\xi_1}}$, for $u\in\ell^2$. Moreover,
$D(A_1)=X=\ell^2=D(A_1)^*$. Obviously $A_1^*=A_1$ and $A_1$ is the generator of a strongly continuous semigroup on $\ell^2$. 

Thus, we may consider $\mathbf{L}$, as a perturbation of $A_1$ by the bounded (skew-adjoint) linear operator $A_2:\ell^2\rightarrow\ell^2$,  $(A_2u)_{n\in\mathbb{Z}}=i\xi_2(Au)_{n\in\mathbb{Z}}$, and apply \cite[Theorem 1.1]{Pazy83}, to justify that $\mathbf{L}=A_1+A_2$, is the generator of a strongly continuous emigroup on $\ell^2$.\ \ $\diamond$  

With Lemma \ref{semgen} at hand, and applying it for the case $\xi_1=\lambda$, $\xi_2=\alpha$, we may recast (\ref{lat1})-(\ref{lat2}) into the form of the integral equation
\begin{eqnarray}
\label{integ1}
u(t)=U(t)u_0+\int_{0}^{t}U(t-s)F_1(u(s))ds,\;\;\; F_1(u):=-(F(u)+\gamma u).
\end{eqnarray}
With the help of Lemma \ref{LipN}, we can handle (\ref{integ1}) by a contraction method  and the local existence result can be stated as follows (we refer to \cite{cazh},\cite{AN} for the proof).
\begin{rigor1}
\label{locex}
We assume that $\lambda,\alpha>0$, and conditions $(\mathrm{N_1})$ or $(\mathrm{N_2})$ are satisfied. Then there exists a function $T^*:\ell^2\rightarrow (0,\infty]$ with the following properties:\vspace{.2cm}\\
(a)\  For all $u_0\in\ell^2$, there exists  $u\in\mathrm{C}([0,T^*(u_0)),\ell^2)$, such that
for all $0<T<T^*(u_0)$, $u$ is the unique solution of (\ref{lat1})-(\ref{lat2}) in $\mathrm{C}([0,T],\ell^2)$ (well posedeness).\vspace{.2cm}\\
(b)\  For all $t\in [0,T^*(u_0))$,
\begin{eqnarray}
\label{maxT}
T^*(u_0)-t\geq \frac{1}{2(L(R)+1)}:=T_R,\;\;R=2||u(t)||_{\ell^2},
\end{eqnarray}
where $L(R)$, is the Lipschitz constant for the map $F_1:\ell^2\rightarrow\ell^2$. Moreover the following alternative holds: (i) $T^{*}(u_0)=\infty$, or (ii) $T^*(u_0)<\infty$ and $\lim_{t\uparrow T^*(u_0)}||u(t)||_{\ell^2}=\infty$ (maximality).\vspace{.2cm}\\
(c)\ $T^*:\ell^2\rightarrow (0,\infty]$ is lower semicontinuous.
In addition, if $\{u_{n0}\}_{n\in\mathbb{N}}$ is a sequence in
$\ell^2$ such that $u_{n0}\rightarrow u_0$ and if $T<T^*(u_0)$,
then $S(t)u_{0n}\rightarrow S(t)u_0$ in $\mathrm{C}([0,T],\ell^2)$, where $S(t)u_0=u(t)$, $t\in [0,T^*(u_0))$, denotes the solution operator
(continuous dependence on initial data).
\end{rigor1}
We note that Lemma \ref{semgen} and Theorem \ref{locex} remains valid in the cases $\lambda=0$, $\alpha >0$ (DNLS type equation) and $\lambda>0$, $\alpha=0$. 
\subsection{Local existence in $\ell^2_{\theta}$ and $\ell^1$}
Since the operator $A$ is not symmetric in  $\ell_\theta^2$, we cannot apply the analysis for the operator
$\mathbf{L}$ and Theorem \ref{locex}, for the local existence in $\ell^2_{\theta}$.
On the other hand for the case of $\ell^1$, Hilbert space methods are not applicable. In both cases, the problem can be treated by general existence Theorems in Banach spaces. 

Concerning the operator $\mathbf{L}:X\rightarrow X$, $X=\ell^2_{\theta},\ell^1$ we have the following 
\begin{smallTheorem}
\label{thm:LOCALLIPS2}
The operator $\mathbf{L} :X\rightarrow X$, defined by
 is globally Lipschitz on $X$.
\end{smallTheorem}
{\bf Proof:}
Let $u,\,v\in B_R$. Then $(\mathbf{L}u)_{n\in\mathbb{Z}}-(\mathbf{L}v)_{n\in\mathbb{Z}}=(\lambda+i\alpha)\left\{(u_{n+1}-v_{n+1})-2(u_n-v_n)+(u_{n-1}-v_{n-1})\right\}$ and it follows that $||Au-Av||_{X}\leq L||u-v||_{X}$, where $L$ is independent of $R$.\ $\diamond$\newline

Thus,the local existence
result in the case of the spaces $X=\ell^2_{\theta},\ell^1$, can be stated as follows.
\begin{rigor1}
\label{lthetaex}
We assume that $\lambda,\alpha>0$, and conditions $(\mathrm{N_1})$ or $(\mathrm{N_2})$ are satisfied. For all $u_0\in X$, there exists
$T^*(u_0)>0$, such that for all $0< T<T^*(u_0)$, there exists a unique solution of the problem (\ref{lat1})-(\ref{lat2}), $u(t)\in \mathrm{C}^1([0,T],X)$.
\end{rigor1}
{\bf Proof:}\ \ This time, we write (\ref{lat1})-(\ref{lat2}), as an ordinary differential equation in $\ell^2_\theta$
\begin{eqnarray*}
\label{absode}
\dot{u}(t)&+&\Phi(u(t))=0,\\
u(0)&=&u_0,\nonumber
\end{eqnarray*}
where $\Phi(u)=F(u)+\mathbf{L}(u)+\gamma u$,  and $u(t)$ lies in
$X$. Lemmas \ref{LipN}-\ref{thm:LOCALLIPS2}
suffice for the application of standard existence and uniqueness
Theorems for ordinary differential equations in Banach spaces (e.g. generalized Peano and Picard-Lindelof Theorems)
\cite[pg. 78-82]{zei85}. \  $\diamond$

We remark that \ref{lthetaex} holds also in the case where $\lambda\geq 0$ and $\alpha,\gamma\in\mathbb{R}$.
\section{The case of a finite lattice with Dirichlet boundary conditions}
We will study in this section the asymptotic behavior of solutions to the DCGL and DRGL equations, considered in a finite lattice, assuming Dirichlet boundary conditions 
\begin{eqnarray}
\label{bvp1}
\dot{u}_n-(\lambda+i\alpha)(u_{n-1}&-&2u_n+u_{n+1})=F(u_n)+\gamma u_n,\;\;|n|\leq N,\\
\label{bvp2}
u_{-(N+1)}(\cdot)&=&u_{(N+1)}(\cdot)=0,\\
\label{bvp3}
u_n(0)&=&u_{n,0},\;\;|n|\leq N.
\end{eqnarray}
Our questions concerning the life span of solutions to (\ref{bvp1})-(\ref{bvp3}) consider nonlinear interactions of the following forms\vspace{.2cm}
\newline
$(\mathcal{NG})$\ \ $F(s)=(k+i\beta)|s|^p$ for some $p>1$. (non-gauge type nonlinearity),\vspace{.2cm}
\newline
$(\mathcal{G})$\ \  $F(s)=(k+i\beta)|s|^{p-1}s$ for some $p>1$ (gauge type nonlinearity). \vspace{.2cm}

First we present theoretical estimates for the blow-up time, for several physically interesting parameter regimes and we conclude with the existence of a global attractor in the case of the gauge nonlinearity.
\subsection{Blow-up in finite time for the case of a DCGL equation in the case of non-gauge nonlinearity}
{\bf A.} ($\lambda\geq 0$ and $\beta>0$).\ \ Motivated by \cite{Ozawa}, for any $t\in (0,T^*)$ we define the function
\begin{eqnarray}
\label{th1}
M(t)=\frac{e^{-\gamma t}}{L^{\sigma}}\mathrm{Im}\sum_{n=-N}^{n=N}u_n(t),\;\; L=2N+1,\;\;\sigma >0,
\end{eqnarray}
and we assume that $\mathrm{Im}\sum_{n=-N}^{n=N}u_{0,n}>0$.
The unspecified parameter $\sigma$, will be related to a scaling argument, on the investigation of the behavior of the upper bound for the blow-up time, that we shall derive in the sequel (see Remark \ref{rem1}). 

We differentiate (\ref{th1}) to obtain
\begin{eqnarray}
\label{th2}
M'(t)&=&\frac{e^{-\gamma t}}{L^{\sigma}}\mathrm{Im}\sum_{n=-N}^{n=N}(-\gamma u_n+\dot{u}_n)=
\frac{e^{-\gamma t}}{L^{\sigma}}\mathrm{Im}\sum_{n=-N}^{n=N}\left\{(\lambda+i\alpha)(u_{n-1}-2u_n+u_{n+1})+F(u_n)\right\}\nonumber\\
&=&\frac{\beta e^{-\gamma t}}{L^{\sigma}}\sum_{n=-N}^{n=N}|u_n|^p\geq 0.
\end{eqnarray}
Now, an application of  inequality (\ref{otherfnorms}) to (\ref{th1}), implies that
\begin{eqnarray*}
M(t)\leq \frac{e^{-\gamma t}}{L^{\sigma}}\mathrm{Im}\sum_{n=-N}^{n=N}|u_n|
\leq\frac{e^{-\gamma t}}{L^{\sigma}}L^{1/q}\left\{\sum_{n=-N}^{n=N}|u_n|^p\right\}^{1/p},
\end{eqnarray*}
with $q=\frac{p}{p-1}$, hence
\begin{eqnarray}
\label{th3}
M(t)^p\leq e^{-p\gamma t}L^{p-1-p\sigma}\sum_{n=-N}^{n=N}|u_n|^p.
\end{eqnarray}
Inserting (\ref{th3}) to (\ref{th2}), we derive the inequality
\begin{eqnarray}
\label{th4}
M'(t)\geq \beta e^{(p-1)\gamma t}L^{(1-p)(1-\sigma)}M(t)^p.
\end{eqnarray}
Now using (\ref{th4}), and differentiating the function $M^{1-p}(t)$ we observe that
\begin{eqnarray}
\label{th5}
\frac{d}{dt}(M(t)^{1-p})\leq -(p-1)\beta e^{(p-1)\gamma t}L^{(1-p)(1-\sigma)}.
\end{eqnarray}
Integration of (\ref{th5}) with respect to time, implies that
\begin{eqnarray}
\label{th6}
M(t)^{1-p}\leq\left\{
\begin{array}{ccc}
M(0)^{1-p}-\frac{\beta}{\gamma}\left(e^{(p-1)\gamma t}-1\right)L^{(1-p)(1-\sigma)},\;\; &\gamma\neq 0,\nonumber \\
M(0)^{1-p}-(p-1)\beta L^{(1-p)(1-\sigma)}t,\;\; &\gamma=0. 
\end{array}
\right.
\end{eqnarray}
Since $M(t)>0$ for all $t\in [0,T^*)$, we obtain from (\ref{th6}), that the maximal existence time $T^*$ can be estimated as
\begin{eqnarray}
\label{th7}
T^*\leq\left\{
\begin{array}{ccc}
\frac{1}{(p-1)\gamma}\ln\left\{1+\frac{\gamma}{\beta}M(0)^{1-p}L^{-(1-p)(1-\sigma)}\right\},\;\; &\gamma\neq 0, \\
M(0)^{1-p}\;\frac{L^{-(1-p)(1-\sigma)}}{(p-1)\beta},\;\; &\gamma=0. 
\end{array}
\right.
\end{eqnarray}
Note that we have assumed that 
\begin{eqnarray}
\label{condiss}
\frac{\gamma}{\beta}M(0)^{1-p}>-L^{(1-p)(1-\sigma)},\;\;\gamma\neq 0. 
\end{eqnarray}
{\bf B.} ($\lambda\geq 0$ and $k>0$)\ \ This time, we consider the quantity
\begin{eqnarray}
\label{kpos}
N(t)=\frac{e^{-\gamma t}}{L^{\sigma}}\mathrm{Re}\sum_{n=-N}^{n=N}u_n(t),\;\; L=2N+1,
\end{eqnarray}
assuming now that $\mathrm{Re}\sum_{n=-N}^{n=N}u_{0,n}>0$.
We observe that 
\begin{eqnarray*}
N'(t)=\frac{k e^{-\gamma t}}{L^{\sigma}}\sum_{n=-N}^{n=N}|u_n|^p\geq 0.
\end{eqnarray*}
Following similar arguments to those we used in case {\bf A.}, we obtain that the maximal existence time $T^*$ can be estimated as
\begin{eqnarray}
\label{th7k}
T^*\leq\left\{
\begin{array}{ccc}
\frac{1}{(p-1)\gamma}\ln\left\{1+\frac{\gamma}{k}N(0)^{1-p}L^{-(1-p)(1-\sigma)}\right\},\;\; &\gamma\neq 0, \\
N(0)^{1-p}\;\frac{L^{-(1-p)(1-\sigma)}}{(p-1)k},\;\; &\gamma=0. 
\end{array}
\right.
\end{eqnarray}
This time we have assumed that 
\begin{eqnarray}
\label{condiss2}
\frac{\gamma}{k}N(0)^{1-p}>-L^{(1-p)(1-\sigma)},\;\;\gamma\neq 0. 
\end{eqnarray}
We summarize the above results, in the following
\begin{rigor1}
\label{acase}
{\bf A.} We assume that $\lambda\geq 0$, $\beta >0$ and 
$\mathrm{Im}\sum_{n=-N}^{n=N}u_{0,n}>0$.
Then for the DCGL equation (\ref{bvp1})-(\ref{bvp3}) with nonlinear interaction $(\mathcal{NG})$, the maximal existence time is estimated by (\ref{th7}).\newline
{\bf B.} We assume that $\lambda\geq 0$, $k>0$ and $\mathrm{Re}\sum_{n=-N}^{n=N}u_{0,n}>0$. Then for the DCGL equation (\ref{bvp1})-(\ref{bvp3}) with nonlinear interaction $(\mathcal{NG})$, the maximal existence time is estimated by (\ref{th7k}).
\end{rigor1}
\begin{remark}
\label{rem1}
{\bf (Scaling limit and Blow-up).} {\em The  estimates for the upper bound of the blow up time $T^{*}$,
can be interpreted in the following way, employing some scaling arguments. The parameter $\sigma$ is unspecified in the above argument. 
Let us consider the behaviour of the upper bound for the blow up time $T^{*}$, as we take the case of a large system ($L\rightarrow \infty$). 
Since $p>1$, and under the assumption that $M(0)=O(1)$,  we see that the upper bound tends to zero if $\sigma >1$, whereas the upper bound tends to infinity if $\sigma <1$. The upper bound is independent of $L$ if $\sigma=1$. This means that if $\sigma >1$, then in the limit of large systems $L\rightarrow \infty$, the system blows up instantly, when $\sigma=1$, one may ask if the blow-up time is independent of the lattice size $L$, whereas if $\sigma <1$, one may ask if the system may have longer lifetimes. 

The above observation, along with the condition that $M(0)=O(1)$, allows for some heuristic investigations, regarding the lifetime of solutions and its possible dependence on the way the initial data decay: from the definition of $M(0)$, we see that this quantity is of order $1$, as long as the sum of $u_n(0)$ on all lattice sites, scales as $L^{\sigma}$. That means that if $u_n(0)\sim n ^{\delta}$ as $n\rightarrow \infty$, then $M(0)\sim L^{1+\delta -\sigma}$, so that $M(0)=O(1)$, as long as $\delta = \sigma -1$. Thus, we may have instant blow up, as long as 
$\delta >0$. On the other hand, it seems to be an interesting question, if the solution may live for longer times as long as $\delta <0$ (since the behavior of the upper bound does not necessarily imply a similar behavior of the  blow-up time). Finally, it appears that the upper bound  is independent of the size of the system, as long as $\delta=0$. In conclusion, if the initial data do not decay fast enough in space, the solution will blow up instantly for large systems, whereas for spatially decaying initial data, one may conjecture that the solution may live for longer finite times.\footnote{However, in practice the actual life time may be smaller, as it can be seen from numerical simulations.}}
\end{remark}
\begin{remark}{\bf (Indications for global existence).} {\em We observe that conditions (\ref{condiss}) and (\ref{condiss2}) are always valid in the case where $\gamma>0$ (i.e in the case where the linear term acts as a linear source). On the other hand in the case $\gamma<0$ (linear dissipation), these conditions (\ref{condiss}) and (\ref{condiss2}) imply that
\begin{eqnarray*}
\label{dissrange}
&&0<-\gamma<\frac{\beta L^{(1-p)(1-\sigma)}}{M(0)^{1-p}},\;\;\beta>0,\\
&&0<-\gamma<\frac{k L^{(1-p)(1-\sigma)}}{N(0)^{1-p}},\;\;k>0,
\end{eqnarray*}
providing a range for the dissipation parameter, for a possible observation of blow-up in finite time.} 
\end{remark}
\subsection{Blow-up in finite time for the case of a DRGL equation in the case of gauge nonlinearity}
We shall examine now the case of a DRGL equation ($\alpha=0$ and $\beta=0$) in the case of the gauge type nonlinearity $(\mathcal{G})$ and $k>0$. For this case, we consider the scalar quantity
\begin{eqnarray}
\label{th8}
E(u)=\frac{1}{L^{\sigma}}\left\{\frac{\lambda}{2}\sum_{n=-N}^{n=N}|(B_du)_n|^2-\frac{\gamma}{2}\sum_{n=-N}^{n=N}|u_n|^2-\frac{k}{p+1}\sum_{n=-N}^{n=N}|u_n|^{p+1}\right\},
\end{eqnarray}
and now we shall consider the function
\begin{eqnarray}
\label{th9}
M(t)=\frac{1}{L^{\sigma}}\sum_{n=-N}^{n=N}|u_n(t)|^2.
\end{eqnarray}
Now we multiply equation (\ref{bvp1}) in the $\mathbb{C}^{2N+1}$ scalar product. We get the energy equation
\begin{eqnarray}
\label{th10}
\frac{1}{2}\frac{d}{dt}\sum_{n=-N}^{n=N}|u_n|^2+\lambda\sum_{n=-N}^{n=N}|(B_du)_n|^2-\gamma\sum_{n=-N}^{n=N}|u_n|^2-k\sum_{n=-N}^{n=N}|u_n|^{p+1}=0.
\end{eqnarray}
Then from (\ref{th9}) and (\ref{th10}) we obtain
\begin{eqnarray}
\label{th11}
M'(t)&=&-\frac{2\lambda}{L^{\sigma}}\sum_{n=-N}^{n=N}|(B_du)_n|^2+\frac{2\gamma}{L^{\sigma}}\sum_{n=-N}^{n=N}|u_n|^2+\frac{2k}{L^{\sigma}}\sum_{n=-N}^{n=N}|u_n|^{p+1}\nonumber\\
&=&-4E(u(t))+\frac{2k(p-1)}{(p+1)L^{\sigma}}\sum_{n=-N}^{n=N}|u_n|^{p+1}.
\end{eqnarray}
Multiplying the DRGL equation by $\frac{1}{L^{\sigma}}\dot{u}_n$ and keeping real parts, we observe that $E(u(t))\leq E(u_0)$. Now under the assumption that $E(u_0)\leq 0$, we get from (\ref{th11})
that
\begin{eqnarray}
\label{th12}
M'(t)\geq -4E(u_0)+\frac{2k(p-1)}{(p+1)L^{\sigma}}\sum_{n=-N}^{n=N}|u_n|^{p+1}\geq \frac{2k(p-1)}{(p+1)L^{\sigma}}\sum_{n=-N}^{n=N}|u_n|^{p+1}.
\end{eqnarray}
Once again, inequality  (\ref{otherfnorms}) implies that
\begin{eqnarray}
\label{finlat}
\sum_{n=-N}^{n=N}|u_n|^2\leq L^{\frac{p-1}{p+1}}\left(\sum_{n=-N}^{n=N}|u_n|^{p+1}\right)^{\frac{2}{p+1}}.
\end{eqnarray}
Therefore, we have for $M(t)$ that
\begin{eqnarray}
\label{th13}
M(t)^{\rho}\leq L^{\frac{p-1-\sigma p-\sigma}{2}}\sum_{n=-N}^{n=N}|u_n|^{p+1},\;\;\rho=\frac{p+1}{2}.
\end{eqnarray}
Now setting $\rho_1=\frac{p+1}{p-1}$, we insert (\ref{th13}) into (\ref{th12}) and we obtain
\begin{eqnarray}
\label{th14}
M'(t)\geq\frac{2k}{\rho_1}L^{(1-p)(1-\sigma)}M(t)^{\rho}.
\end{eqnarray}
We proceed as for the derivation of the estimate (\ref{th7}). Since $\rho>1$, we have from (\ref{th14}) that
\begin{eqnarray*}
\frac{d}{dt}(M(t)^{1-\rho})\leq -(\rho -1)\frac{2k}{\rho_1}L^{(1-p)(1-\sigma)},
\end{eqnarray*}
and integration with respect to time, implies that
\begin{eqnarray*}
M(t)^{1-\rho}\leq M(0)^{1-\rho}-(\rho-1))\frac{2k}{\rho_1}L^{(1-p)(1-\sigma)}t.
\end{eqnarray*}
Since $M(t)\geq 0$ for all $t\in [0,T^*)$, we have that this time, the maximal existence time $T^*$ satisfies the estimate
\begin{eqnarray}
\label{th15}
T^*\leq \frac{(p+1)}{k(p-1)^2}\frac{L^{-(1-p)(1-\sigma)}}{M(0)^{\frac{p-1}{2}}}.
\end{eqnarray}
Summarizing, in the case of the gauge nonlinearity, we have the following 
\begin{rigor1}
\label{bcase}
We assume that $\lambda>0$, $k >0$ and that the initial energy is such that $E(u_0)\leq 0$. 
Then for the DRGL equation (\ref{bvp1})-(\ref{bvp3}) ($\alpha=\beta=0$) with nonlinear interaction $(\mathcal{G})$, the maximal existence time is estimated by (\ref{th15}).
\end{rigor1}
\begin{remark}(Scaling limit and Blow-up) {\em Similar comments as in Remark \ref{rem1} hold with respect to scaling of the initial data and blow-up, for the case of the DRGL equation (\ref{bvp1})-(\ref{bvp3}), with the gauge nonlinear interaction $(\mathcal{G})$.}
\end{remark}
\subsection{Comparison of  DRGL equations with DNLS, with respect to global existence of solutions and blow-up in finite time}\label{sec:MAGIRIO}
The fact that solutions of DRGL, in the case of gauge nonlinearity $(\mathcal{G})$, may blow-up in finite time, is in  contrast with the behavior of solutions of the discrete nonlinear Schr\"odinger equation (DNLS) with the same nonlinearity. The solutions of DNLS in the case of gauge nonlinearity, exist globally, unconditionally with respect to the degree of the nonlinearity, the size of the initial data and the sign of the initial energy as it was observed first numerically in \cite{bang}. 
A detailed discussion on the asymptotic behavior of solutions of DNLS equations is presented\cite{AN}. For the sake of completeness and for a comparison, we present here the simple proof on global existence of solutions, in the case of the infinite lattice. That is, we shall consider the DNLS lattice differential equation (the case $\lambda=k=0$ of (\ref{lat1})-(\ref{lat2}))
\begin{eqnarray}
\label{cons1}
\dot{u}_n-i\alpha(u_{n-1}&-&2u_n+u_{n+1}) =i\beta|u_n|^{p-1}u_n,\;\;n\in\mathbb{Z},\;\;1< p <\infty,\\
\label{cons2}
u_n(0)&=&u_{n,0},\;\;n\in\mathbb{Z}.
\end{eqnarray}   
and we assume that $\alpha,\beta>0$. Theorem \ref{locex} (or \ref{lthetaex}) covers also the case of DNLS (\ref{cons1})-(\ref{cons2}): For all $u_0\in\ell^2$, there exists  $u\in\mathrm{C}([0,T^*(u_0)),\ell^2)$ such that
for all $0<T<T^*(u_0)$, $u$ is the unique solution of (\ref{cons1})-(\ref{cons2}) in $\mathrm{C}([0,T],\ell^2)$. 
Taking the scalar product of (\ref{cons1}) with $iu$, we obtain that
\begin{eqnarray}
\label{charge}
\frac{d}{dt}||u(t)||^2_{\ell^2}=0,\;\;\mbox{or}\;\;||u(t)||^2_{\ell^2}=||u_0||^2_{\ell^2},\;\;\mbox{for every}\;\;t\in [0,T^*(u_0)),
\end{eqnarray}
Although the conserved quantity (\ref{charge}) suffices to demonstrate global existence, to elucidate the interplay of nonlinearity and discreteness, we shall examine the DNLS Hamiltonian 
\begin{eqnarray}
\label{eneg2a}
\mathrm{E}(u(t))=\mathrm{E}(u_0),\;\;\mathrm{E}(u(t)):=\frac{\alpha}{2}\sum_{n\in\mathbb{Z}}|(Bu)_n(t)|^2-
\frac{\beta}{p+1}\sum_{n\in\mathbb{Z}}|u_n(t)|^{p+1}.
\end{eqnarray}
From (\ref{charge}) and (\ref{eneg2a}) we may easily derive the ''conservation law" 
\begin{eqnarray}
\label{eneg2}
\mathrm{E_1}(u(t))=\mathrm{E}_1(u_0),\;\;\mathrm{E}_1(u(t)):=\frac{\alpha}{2}||u(t)||_{\ell^2_1}^2-
\frac{\beta}{p+1}||u(t)||_{\ell^{p+1}}^{p+1},
\end{eqnarray}
where the $\ell^2_1$ is defined by (\ref{lp4}), for $\xi_1=1$. 
Then, by using  (\ref{lp1}), (\ref{charge}) and (\ref{eneg2}), we may derive the estimate
\begin{eqnarray}
\label{lastest}
||u(t)||_{\ell^{2}_1}^2
&\leq& ||u_0||_{\ell^2_1}^2+\frac{2\beta}{\alpha(p +1)}\left\{||u_0||^{p+1}_{\ell^{p+1}}+||u(t)||^{p+1}_{\ell^{p+1}}\right\}\nonumber\\
&\leq& ||u_0||_{\ell^2_1}^2+\frac{2\beta}{\alpha(p+1)}\left\{||u_0||^{p+1}_{\ell^2}+||u(t)||^{p+1}_{\ell^2}\right\}\nonumber\\
&\leq & ||u_0||_{\ell^2_1}^2+ \frac{4\beta}{\alpha(p+1)}||u_0||^{p+1}_{\ell^2}.
\end{eqnarray}
As a consequence of (\ref{lastest}) we obtain  that $T^*(u_0)=\infty$ and $\mathrm{sup}\left\{||u(t)||_{\ell^2_1},\,t\in [0,\infty)\right\}<\infty$. The proof is very similar in the case of the finite lattice, assuming Dirichlet boundary conditions.

Let us mention that this behavior of the DLNS system, is not only in contrast with the DRGL system,  but also with its continuous counterpart, 
\begin{eqnarray}
\label{pde}
\partial_t u -i\alpha u_{xx}&=&i\beta|u|^{p-1}u,\;\; x\in\mathbb{R},\;\;t>0,\\
u(x,0)&=&u_0(x).\nonumber
\end{eqnarray}
In order to clarify these differences, let us recall  the main results concerning (\ref{pde})
(see \cite{cazh,cazS,Yvan}): For $u_0\in\mathrm{H^1}(\mathbb{R})$ and $1 <p <\infty$ there exists a unique maximal solution of (\ref{pde}), $u(t)\in \mathrm{C}([0,T_{max}),\mathrm{H^1}(\mathbb{R}))\cap\mathrm{C}^1([0,T_{max}),\mathrm{L^2}(\mathbb{R}))$. In addition:
If $1< p<5$ then $T_{\max}=\infty$ and $u$ is bounded in $\mathrm{H^1}(\mathbb{R})$.
Let $p\geq 5$. Assume that $u_0\in \mathrm{H^1}(\mathbb{R})$ such that $\int_{\mathbb{R}}|x|^2|u_0|^2dx <\infty$ (initial data with finite variance) and $E(u_0)<\infty$. Then $T_{\max}<\infty$.
On the other hand if $||u_0||_{\mathrm{H^1}}$ is sufficiently small, $T_{max}=\infty$ and $u$ is bounded in $\mathrm{H^1}(\mathbb{R})$.

In the case of DNLS, the assumption of ininitial data with finite variance, reads as $\sum_{n\in\mathbb{Z}}|n|^2|u_{n,0}|^2<\infty$, and such data belong to $\ell^2$. On the other hand, it follows from Theorem \ref{acase} A., that the solutions of DNLS system with non-gauge nonlinear interaction $(\mathcal{NG})$ may blow up in finite time (at least in the case of Dirichlet boundary conditions). Thus, regarding the DNLS system, with respect to global existence of solutions and blow-up in finite time we may comment with the following
\begin{figure}
\centerline{\includegraphics[angle=0,width=0.75\textwidth]{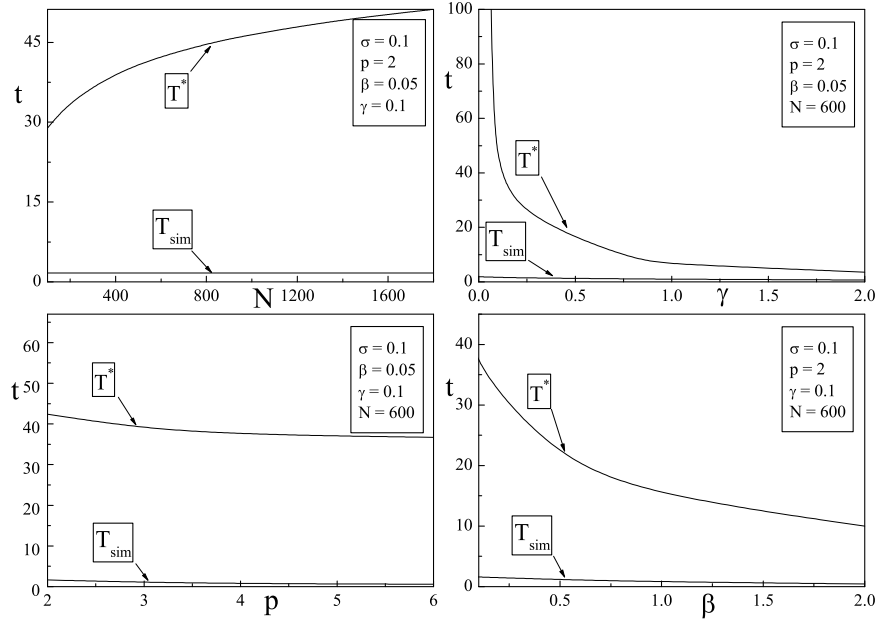}}
\caption{Theoretical upper bound for the blow-up time ($T^{*}$) and {\em numerically obtained ($T_{\mathrm{sim}}$) blow-up time} for the finite lattice, in the case of initial data for which $\sigma < 1$.}
\label{fig:LATTICE1}
\end{figure}
\begin{rigor1}
\label{DNLSBL} 
{\bf A.} Consider the DNLS equation (\ref{cons1}) with Dirichlet boundary conditions and nonlinear interaction $(\mathcal{NG})$.  Assume that  $\alpha ,\beta >0$ and $\mathrm{Im}\sum_{n=-N}^{n=N}u_{0,n}>0$.
Then the solution blows-up in finite time for all $1<p<\infty$. The maximal existence time is estimated by (\ref{th7}).\newline
{\bf B.} Consider the DNLS equation (\ref{cons1}) with Dirichlet boundary conditions and nonlinear interaction $(\mathcal{G})$. Assume that  $\alpha ,\beta >0$. Then the solution exists globally in time unconditionally with respect to the initial data and the sign of the initial energy, for all $1<p<\infty$. 
The  same holds for the DNLS infinite lattice That is for any $u_0\in\ell^2$, the solution of (\ref{cons1}) is in $\mathrm{C}^1([0,\infty),\ell^2)$.
\end{rigor1}

\subsection{Numerical simulations}
Eventhough the aim of the paper, was to show -using some analytical arguments-finite time blow-up of solutions of discrete complex lattices, we decided to test the
theoretical estimates for the upper bound of the blow-up time, and the heuristic scaling arguments of Remark \ref{rem1}, numerically, against the observed blow-up times, for some parameter values. 

The finite lattice equations, have been integrated numerically over time, using a fourth order Runge Kutta scheme, implementing Dirichlet boundary conditions. For the case of non-gauge type nonlinearity, and for initial data with $\sigma<1$, the {\em numerically obtained blow up time $T_{\mathrm{sim}}$}, is shown in figure  \ref{fig:LATTICE1} and compared with the theoretical estimate for the upper bound of $T^{*}$ (still denoted for simplicity, by $T^{*}$). 

In the first graph, we show the variation of $T_{\mathrm{sim}}$ and the upper bound for $T^{*}$,  as a function of the number of lattice sites $N$, in the second with respect to $\gamma$, in the third with respect to $p$ and in the fourth with respect to $\beta$. 
In the first graph, it seems that the observed blow-up time is independent of $N$ (although inspection of the data, show a slow increase).  Note however, that according to (\ref{th7}), and  as clearly stated in Remark \ref{rem1}, when $\sigma<1$, the behavior of the upper bound, does not necessarily imply a similar behavior of the blow-up time.
\begin{figure}
\centerline{\includegraphics[angle=0,width=0.75\textwidth]{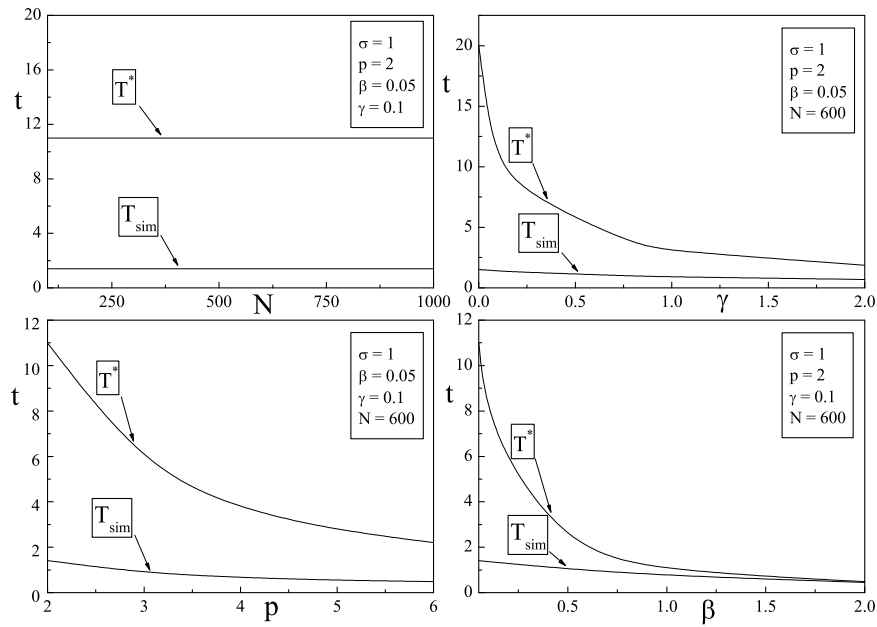}}
\caption{Theoretical upper bound for the blow-up time ($T^{*}$), and {\em numerically obtained ($T_{\mathrm{sim}}$) 
blow-up time} for the finite lattice in the case of initial data for which $\sigma = 1$.}
\label{fig:LATTICE2}
\end{figure}
In figure \ref{fig:LATTICE2}, we demonstrate the case of initial data with $\sigma=1$. We first observe, that the scaling limit argument we propose, with respect to the number of lattice sites, and  the theoretical estimates for the variation of the upper bound ,with respect to the various parameters discussed above, seem to capture-at least quatitatively- the variation of the numerically observed blow-up times, with respect to these parameters: in the case $\sigma=1$, the numerically obtained blow-up time seem to be independent of the lattice size $L=2N+1$.  
\begin{figure}
\centerline{\includegraphics[angle=0,width=0.75\textwidth]{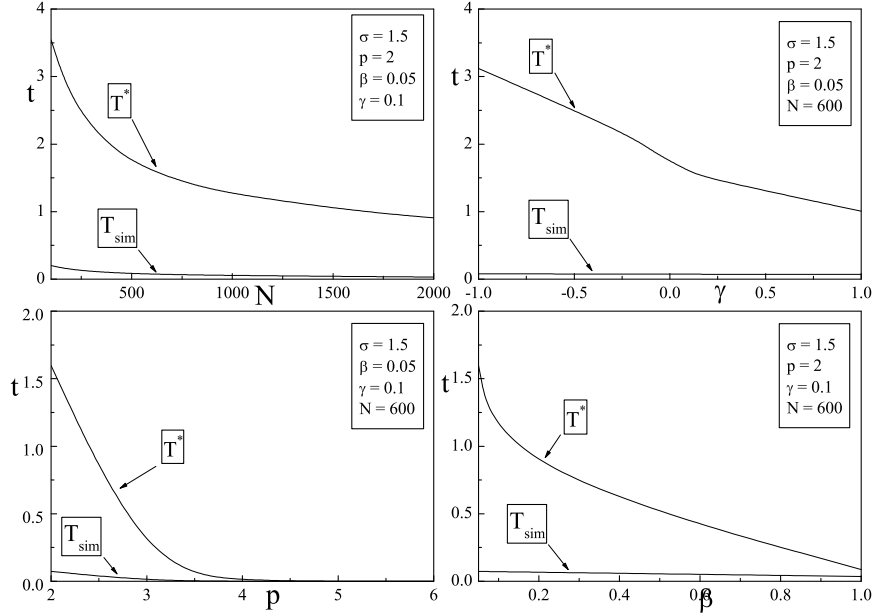}}
\caption{Theoretical upper bound for the blow-up time  ($T^{*}$), and {\em numerically obtained ($T_{\mathrm{sim}}$) blow-up time} for the finite lattice in the case of initial data for which $\sigma > 1$.}
\label{fig:LATTICE3}
\end{figure}
Furthermore, regarding the variation with respect to the parameters (and especially with parameter $\beta$), it is surprising that the upper bound is close (same order of magnitude) with the numerically observed blow-up times.  We also observe that the trend predicted by the theory seem to be verified, as far as the variation with respect to the parameters is concerned.

Finally in figure \ref{fig:LATTICE3} we repeat the same calculations with initial data with $\sigma>1$. Especially in the case of the variation with respect to the number of the lattice sites, the numerical simulations seem to verify the instant blow-up of solutions for increasing lattice size, as it is predicted by the proposed scaling argument. 
\subsection{A global attractor for the case of non-gauge nonlinearity}
We conclude our study for the finite complex lattice, by verifying existence of a global attractor, for the case of a dissipative gauge nonlinearity.
\begin{rigor2}
\label{fattr}
We assume that $\lambda>0$, $\alpha,\beta,\gamma\in\mathbb{R}$ and $k<0$. Let $u_0:=(u_{n,0})_{|n|\leq N}\in\mathbb{C}^{2N+1}$. For (\ref{bvp1})-(\ref{bvp3}), with nonlinear interaction $(\mathcal{G})$, there exists a unique solution
(\ref{bvp1})-(\ref{bvp3})
such that $u\in\mathrm{C}^1([0,\infty),\mathbb{C}^{2N+1})$. The dynamical system defined by  (\ref{bvp1})-(\ref{bvp3})
\begin{eqnarray}
\label{fdynamical} \mathcal{S}_N(t):u_0\in \mathbb{C}^{2N+1}\rightarrow
u(t)\in \mathbb{C}^{2N+1},
\end{eqnarray}
possesses a bounded absorbing set $\mathcal{O}_0$ in
$\mathbb{C}^{2N+1}$ and a global attractor
$\mathcal{A}_{N}=\omega (\mathcal{O}_0) \subset\mathcal{O}_0\subset\mathbb{C}^{2N+1}$. Moreover, for the absorbing ball, the entry time is independent of the initial data in $\mathbb{C}^{2N+1}$ , however large is the norm of the initial data. 
\end{rigor2}
{\it Proof:}\ \ Let $k=-m,\;m>0$. Taking the scalar product of (\ref{bvp1}) with $u$, we obtain the equation
\begin{eqnarray}
\label{finlatenerg}
\frac{1}{2}\frac{d}{dt}\sum_{n=-N}^{n=N}|u_n|^{2}+\lambda\sum_{n=-N}^{n=N}|(B_du)_n|^2-\gamma\sum_{n=-N}^{n=N}|u_n|^{2} +m\sum_{n=-N}^{n=N}|u_n|^{p+1}=0.
\end{eqnarray}
The operator $-A_d:(\mathbb{C}^{2N+1},||\cdot||_2)\rightarrow (\mathbb{C}^{2N+1},||\cdot||_2)$, defined by (\ref{finself}), is self-adjoint.  If $\lambda_1^*>0$ denotes the smallest eigenvalue for the eigenvalue problem $-A_d \psi=\lambda^*\psi$, then since
$$\lambda_1^*=\mathrm{\inf}_{\psi\in\mathbb{C}^{2N+1}}\frac{(-A_d\psi,\psi)_2}{(\psi,\psi)_2}=
\mathrm{\inf}_{\phi\in\mathbb{C}^{2N+1}}\frac{\sum_{n=-N}^{n=N}|\psi_{n+1}-\psi_n|^2}{\sum_{n=-N}^{n=N}|\psi_n|^2},$$
it is straightforward to check that if $\frac{\gamma}{\lambda}<\lambda_1^*$, the dynamics are trivial in the sense $\limsup_{t\rightarrow\infty}||u(t)||_2=0$. Thus we consider the case $\frac{\gamma}{\lambda}>\lambda_1^*$.
From inequality (\ref{finlat}) we get that
\begin{eqnarray*}
k_1||u||_2^{p+1}\leq m\sum_{n=-N}^{n=N}|u_n|^{p+1},\;\;k_1=m{L}^{\frac{1-p}{2}},
\end{eqnarray*}
while from Young's inequality we get that
\begin{eqnarray*}
\gamma||u||_2^2\leq \frac{k_1}{2}||u||_{2}^{p+1}+\rho_0,\;\;\rho_0(N,\gamma,p)=\frac{p-1}{p+1}\left(\frac{2}{k_1(p+1)}\right)^{\frac{1}{p-1}}\gamma^{\frac{p+1}{p-1}}.
\end{eqnarray*}
Thus, (\ref{finlatenerg}) becomes
\begin{eqnarray}
\label{pregronw}
\frac{1}{2}\frac{d}{dt}||u||_{2}^2+\frac{k_1}{2}||u||_{2}^{p+1}\leq \rho_0.
\end{eqnarray}
Now the result follows by applying Gronwall's Lemma \cite[Lemma 5.1, pg 167-168]{RTem88}: We get from
(\ref{pregronw}) that
\begin{eqnarray}
\label{gronw1}
||u(t)||_2^2\leq\left(\frac{2\rho_0}{k_1}\right)^{\frac{1}{p+1}}+\left(\frac{2}{k_1pt}\right)^{\frac{1}{p}},\;\;\mbox{for all}\;\;t>0.
\end{eqnarray}
Now for any $\rho_1$ satisfying 
$$\rho_1>\left(\frac{2\rho_0}{k_1}\right)^{\frac{1}{p+1}},$$ we derive from (\ref{gronw1}), that
for any set $\mathcal{O}$ of $\mathbb{C}^{2N+1}$, it holds $\mathcal{S}_N(t)\mathcal{O}\subset\mathcal{O}_0$ for any
$$t\geq t_0:=\frac{2}{k_1p}\left(\rho_1^2-\left(\frac{2\rho_0}{k_1}\right)^{\frac{1}{p+1}}\right)^{-p},$$
where $\mathcal{O}_0$ is the ball of $\mathbb{C}^{2N+1}$ of center $0$ and radius $\rho_1$. \ \ $\diamond$.
\begin{remark}
{\em In the case of a DRGL equation with  nonlinearity $(\mathcal{G})$ and $k<0$, the functional $E$ defined by (\ref{th8}) is a Lyapunov function. Moreover, it can be checked (by following similar calculations as those of Proposition \ref{fattr}), that the set of stationary points $\mathcal{E}$ is bounded. Hence, as it follows from \cite{Ball2,jhale88,RTem88},
for each complete orbit  containing $u_0$ lying in
$\mathcal{A}_N$, the limit set $\omega(u_0)$ is a
connected subset of $\mathcal{E}$, on which $E$ is
constant. If $\mathcal{E}$ is totally disconnected (in particular
if $\mathcal{E}$ is countable), any solution $\mathcal{S}_N(t)u_0$ tends to an equilibrium point as $t\rightarrow+\infty$. However, even  in this finite dimensional problem for the DRGL, it could be an interesting issue, the investigation and rigorus justification of the bifuractions from the eigenvalue $\lambda^*_1$, and convergence to (nontrivial) equilibria: writing the stationary DRGL problem as a real system, consisting of the (nonlinear algebraic) equations for the real part $\mathrm{Re}(u)$ and imaginary part $\mathrm{Im}(u)$, one could ask if the eigenvalue $\lambda_1^*$ could be a
bifurcation point, from which two global branches bifurcate. These
branches could consist of semitrivial solutions (i.e. solutions of the form
$(\mathrm{Re}(u),0)$ or $(0,\mathrm{Im}(u))$. Questions of this type will be considered elsewhere \cite{ThNi56}.}
\end{remark} 
\section{The case of an infinite lattice: Existence of global attractor for exponentially spatially localized solutions}
In this section, we prove the existence of a global attractor for the following complex lattice differential equation
\begin{eqnarray}
\label{latg1}
i\dot{u}_n+(\hat{\alpha}+i\hat{\beta})(u_{n-1}&-&2u_n+u_{n+1})+(\hat{\gamma}+i\hat{\delta})u_n+(\hat{\eta}+i\hat{\zeta})F(u_n)=g_n,\;\;n\in\mathbb{Z},\\
\label{latg2}
u_n(0)&=&u_{n,0},\;\;n\in\mathbb{Z}.
\end{eqnarray}
We focus on the case of a gauge nonlinear interaction $F(s)=|s|^{p-1}s,\;p>1$. For specific values of the parameters, one recovers either the DCGL and DRGL equation or the  weakly damped and driven DNLS. We refer to the pioneering work \cite{Bates}, on the existence of global attractors for lattice dynamical systems of first order and in \cite{SZ2} for extensions to various multidimensional lattices of first and second order. We remark that the discretization of the Laplacian, {\em is not self-adjoint in $\ell^2_{\theta}$}, a difference with \cite{Bates} and the examples provided in \cite{SZ2}, which increases considerably the manipulations needed, for the derivation of suitable estimates. The first result, is for the existence of an absorbing ball.  
\begin{smallTheorem}
\label{ballweighted}
Assume condition $\mathrm{(WS)}$ on the weight function and that the parameters $\hat{\alpha},\hat{\beta},\hat{\gamma},\hat{\delta},\hat{\eta},\hat{\zeta}$ are chosen such that for some fixed $\epsilon>0$,
\begin{eqnarray}
\label{parcrit0}
\sigma_0:=\hat{\delta}-\frac{\epsilon}{2}-2\hat{\beta}-|\hat{\alpha}|D\underline{d}^{-1/2}-|\hat{\beta}|\left(1+\frac{\overline{d}}{2}+\frac{\underline{d}^{-1}}{2}\right)>0,\;\;\hat{\zeta}>0.
\end{eqnarray}
Let $(u_{0,n})_{n\in\mathbb{Z}}=u_0\in\ell^2_{\theta}$ and $(g_n)_{n\in\mathbb{Z}}=g\in\ell^2_{\theta}$. A dynamical system can be defined by
(\ref{latg1})-(\ref{latg2}),
\begin{eqnarray}
\label{dynamicalw} S(t):u_0\in {\ell}^2_{\theta}\rightarrow u(t)\in
{\ell}^2_{\theta},
\end{eqnarray}
possessing a bounded absorbing set $\mathcal{B}_0$ in ${\ell}^2_{\theta}$:
For every bounded set $\mathcal{B}$ of ${\ell}^2_{\theta}$, there exists
$t_0(\mathcal{B},\mathcal{B}_0)$ such that for all  $t\geq
t_0(\mathcal{B},\mathcal{B}_0)$, it holds
$S(t)\mathcal{B}\subset\mathcal{B}_0$.
\end{smallTheorem}
{\bf Proof:} We multiply (\ref{lat1}) with $\theta_n\overline{u}_n$, $n\in\mathbb{Z}$ add over all
lattice sites, and keep the imaginary  part. We obtain the equation
\begin{eqnarray}
\label{attract1}
\frac{1}{2}\frac{d}{dt}\sum_{n\in\mathbb{Z}}\theta_n|u_n|^2-\hat{\alpha}\mathrm{I}_1(u_n)-\hat{\beta} \mathrm{I}_2(u_n)+\hat{\delta}\sum_{n\in\mathbb{Z}}\theta_n|u_n|^2
+\hat{\zeta}\sum_{n\in\mathbb{Z}}\theta_n|u_n|^{p+1}=\mathrm{Im}\sum_{n\in\mathbb{Z}}\theta_ng_n\overline{u}_n,
\end{eqnarray}
where the terms $\mathrm{I}_1,\mathrm{I}_2$, are defined as
\begin{eqnarray*}
\label{attract2}
\mathrm{I}_1(u_n)=\sum_{n\in\mathbb{Z}}\left\{(Bu_1)_nB(\theta u_1)_n+(Bu_2)_n(B\theta u_2)_n\right\}
&=&\sum_{n\in\mathbb{Z}}(\theta_{n+1}-\theta_n)(u_{1,n}u_{2,n+1}-u_{2,n}u_{1,n+1}),\;\;\;\;\;\;\;\;\\
\label{attract3}
\mathrm{I}_2(u_n)=\sum_{n\in\mathbb{Z}}\left\{(Bu_1)_n(B\theta u_1)_n+(Bu_2)_n(B\theta u_2)_n\right\}
&=&2\sum_{n\in\mathbb{Z}}\theta_n(u_{1,n}^2+u_{2,n}^2)\\
&&-\sum_{n\in\mathbb{Z}}(\theta_{n+1}+\theta_n)(u_{1,n}u_{1,n+1}+u_{2,n}u_{2,n+1}).
\end{eqnarray*}
Using $\mathrm{(WS)}$, for the term $\mathrm{I}_1$, we may get the estimate
\begin{eqnarray}
\label{attract4}
|I_1(u_n)|&\leq& \sum_{n\in\mathbb{Z}}|\theta_{n+1}-\theta_{n}|\,|u_{1,n}u_{2,n+1}-u_{2,n}u_{1,n+1}|
\leq D\sum_{n\in\mathbb{Z}}\theta_n|u_{1,n}u_{2,n+1}-u_{2,n}u_{1,n+1}|\nonumber\\
&\leq& D\left\{\sum_{n\in\mathbb{Z}}\theta_n|u_{1,n}u_{2,n+1}|+ \sum_{n\in\mathbb{Z}}\theta_n|u_{2,n}u_{1,n+1}|\right\}\nonumber\\
&\leq&
D\left\{\left(\sum_{n\in\mathbb{Z}}\theta_n|u_{1,n}|^2\right)^{1/2}
\left(\sum_{n\in\mathbb{Z}}\theta_n|u_{2,n+1}|^2\right)^{1/2}
+\left(\sum_{n\in\mathbb{Z}}\theta_n|u_{2,n}|^2\right)^{1/2}
\left(\sum_{n\in\mathbb{Z}}\theta_n|u_{1,n+1}|^2\right)^{1/2}\right\}\nonumber\\
&\leq&
D\left\{\left(\sum_{n\in\mathbb{Z}}\theta_n|u_{1,n}|^2\right)^{1/2}
\underline{d}^{-1/2}\left(\sum_{n\in\mathbb{Z}}\theta_n|u_{2,n}|^2\right)^{1/2}
+\left(\sum_{n\in\mathbb{Z}}\theta_n|u_{2,n}|^2\right)^{1/2}
\underline{d}^{-1/2}\left(\sum_{n\in\mathbb{Z}}\theta_n|u_{1,n}|^2\right)^{1/2}\right\}\nonumber\\
&\leq&
2D\underline{d}^{-1/2}\left(\sum_{n\in\mathbb{Z}}\theta_n|u_{1,n}|^2\right)^{1/2}
\left(\sum_{n\in\mathbb{Z}}\theta_n|u_{2,n}|^2\right)^{1/2}\nonumber\\
&\leq& D\underline{d}^{-1/2}||u||_{\ell^2_{\theta}}^2.
\end{eqnarray}
For the second term on the rhs of $\mathrm{I}_2$, we have
\begin{eqnarray}
\label{attract5}
&&\left|\sum_{n\in\mathbb{Z}}(\theta_{n+1}+\theta_{n})(u_{1,n}u_{1,n+1}+u_{2,n}u_{2,n+1})\right|
\leq\sum_{n\in\mathbb{Z}}\theta_{n+1}|u_{1,n}u_{1,n+1}|+\sum_{n\in\mathbb{Z}}\theta_{n}|u_{1,n}u_{1,n+1}|\nonumber\\
&&+\sum_{n\in\mathbb{Z}}\theta_{n+1}|u_{2,n}u_{2,n+1}|+\sum_{n\in\mathbb{Z}}\theta_{n}|u_{2,n}u_{2,n+1}|\nonumber\\
&\leq&
\left(\sum_{n\in\mathbb{Z}}\theta_{n+1}|u_{1,n}|^2\right)^{1/2}
\left(\sum_{n\in\mathbb{Z}}\theta_{n+1}|u_{1,n+1}|^2\right)^{1/2}
+
\left(\sum_{n\in\mathbb{Z}}\theta_{n}|u_{1,n}|^2\right)^{1/2}
\left(\sum_{n\in\mathbb{Z}}\theta_{n}|u_{1,n+1}|^2\right)^{1/2}\nonumber\\
&&+
\left(\sum_{n\in\mathbb{Z}}\theta_{n+1}|u_{2,n}|^2\right)^{1/2}
\left(\sum_{n\in\mathbb{Z}}\theta_{n+1}|u_{2,n+1}|^2\right)^{1/2}
+
\left(\sum_{n\in\mathbb{Z}}\theta_{n}|u_{2,n}|^2\right)^{1/2}
\left(\sum_{n\in\mathbb{Z}}\theta_{n}|u_{2,n+1}|^2\right)^{1/2}\nonumber\\
&\leq&
\frac{1}{2}\left\{\sum_{n\in\mathbb{Z}}\theta_{n+1}|u_{1,n}|^2+\sum_{n\in\mathbb{Z}}\theta_{n+1}|u_{1,n+1}|^2
+\sum_{n\in\mathbb{Z}}\theta_{n}|u_{1,n}|^2+\sum_{n\in\mathbb{Z}}\theta_{n}|u_{1,n+1}|^2\right\}\nonumber\\
&&+\frac{1}{2}\left\{\sum_{n\in\mathbb{Z}}\theta_{n+1}|u_{2,n}|^2+\sum_{n\in\mathbb{Z}}\theta_{n+1}|u_{2,n+1}|^2
+\sum_{n\in\mathbb{Z}}\theta_{n}|u_{2,n}|^2+\sum_{n\in\mathbb{Z}}\theta_{n}|u_{2,n+1}|^2\right\}\nonumber\\
&\leq&
\frac{1}{2}\left\{\overline{d}\sum_{n\in\mathbb{Z}}\theta_{n}|u_{1,n}|^2+\sum_{n\in\mathbb{Z}}\theta_{n+1}|u_{1,n+1}|^2
+\sum_{n\in\mathbb{Z}}\theta_{n}|u_{1,n}|^2+\underline{d}^{-1}\sum_{n\in\mathbb{Z}}\theta_{n+1}|u_{1,n+1}|^2\right\}\nonumber\\
&&+\frac{1}{2}\left\{\overline{d}\sum_{n\in\mathbb{Z}}\theta_{n}|u_{2,n}|^2+\sum_{n\in\mathbb{Z}}\theta_{n+1}|u_{2,n+1}|^2
+\sum_{n\in\mathbb{Z}}\theta_{n}|u_{2,n}|^2+\underline{d}^{-1}\sum_{n\in\mathbb{Z}}\theta_{n+1}|u_{2,n+1}|^2\right\}\nonumber\\
&=&\left(1+\frac{\overline{d}}{2}+\frac{\underline{d}^{-1}}{2}\right)||u||_{\ell^2_{\theta}}^2.
\end{eqnarray}
We insert (\ref{attract4}), (\ref{attract5}) to (\ref{attract1}), to get the inequality
\begin{eqnarray}
\label{absg1}
\frac{d}{dt}||u||^2_{\ell^2_{\theta}}+2\sigma_0||u||^2_{\ell^2_{\theta}}\leq \frac{1}{\epsilon}||g||_{\ell^2_{\theta}}^2.
\end{eqnarray}
From (\ref{absg1}) we derive that
$u\in\mathrm{L^{\infty}}([0,\infty),\ell^2_{\theta})$: Gronwall's Lemma
implies that
\begin{eqnarray}
\label{gronw}
||u(t)||_{\ell^2_{\theta}}^2\leq ||u_0||^2_{\ell^2_{\theta}}\exp(-2\sigma_0 t)+\frac{1}{2\sigma_0\epsilon}||g||_{\ell^2_{\theta}}^2
\{1-\exp(-2\sigma_0 t)\}.
\end{eqnarray}
Letting $t\rightarrow\infty$ we infer that
\begin{eqnarray*}
\label{abbset2}
\limsup_{t\rightarrow\infty}||u(t)||^2_{\ell^2_{\theta}}\leq\frac{1}{2\sigma_0\epsilon}||g||_{\ell^2_{\theta}}^2.
\end{eqnarray*}
Setting $\rho^2=||g||^2_{\ell^2_{\theta}}/2\sigma_0\epsilon$, it follows that for any number $\rho_1>\rho$ the ball $\mathcal{B}_0$ of $\ell^2_{\theta}$ centered at $0$ of radius $\rho_1$ is an absorbing set
for the semigroup $S(t)$: if $\mathcal{B}$ is a bounded set of $\ell^2_{\theta}$, included in a ball of $\ell^2_{\theta}$ centered at $0$ of radius $R$, then for 
$t\geq t_0(\mathcal{B},\mathcal{B}_0)$ where
\begin{eqnarray}
\label{time}
t_0=\frac{1}{2\sigma_0}\log\frac{R^2}{\rho_1^2-\rho^2},
\end{eqnarray}
it holds $||u(t)||^2_{\ell^2_{\theta}}\leq\rho_1^2$, i.e. $S(t)\mathcal{B}\subset\mathcal{B}_0$.
Note that in the absence of external excitation, the dynamical system exhibits trivial dynamics, in the sense that  $\limsup_{t\rightarrow\infty}||u(t)||_{\ell^2_{\theta}}^2=0$, as it follows from (\ref{gronw}).\ \ $\diamond$

The next lemma provides us with the appropriate estimates, on the tail ends of solutions of ({\ref{latg1})-(\ref{latg2}). 
\begin{smallTheorem}
\label{DCGLtail}
Let $u_0\in\mathcal{B}$ where $\mathcal{B}$ is a bounded set of $\ell^2_{\theta}$, and $g\in\ell^2_{\theta}$. Moreover, we assume that the condition (\ref{parcrit0}) on the parameters, is satisfied.
Then, for any $\eta >0$, there exist $T(\eta )$ and $K(\eta )$ such that the solution
$u$ of (\ref{latg1})-(\ref{latg2}) satisfies  for all $t\geq T(\eta )$, the estimate
\begin{eqnarray}
\label{precom}
\sum_{\mid n\mid >2M}\theta_n\mid u_n(t)\mid^2 \le \frac{\eta}{\sigma_0},\;\;\mbox{for any}\;\;
M>K(\eta ).
\end{eqnarray}
\end{smallTheorem}
{\bf Proof:}\ 
We consider a smooth function $\phi \in C^{1}({\mathbb R}^{+},{\mathbb R})$, satisfying the following properties
\begin{eqnarray}
\left\{
\begin{array}{ccc}
\phi(s)=0,\;\; & 0\le s\le 1  \nonumber \\
0\le \phi(s) \le 1,\;\; & 1\le s\le 2 \nonumber \\
\phi(s)=1,\;\; & s\ge 2. \nonumber
\end{array}
\right.
\end{eqnarray}
and
\begin{eqnarray}
\label{mvth}
\mid \phi^{'}(s)\mid \le C_0, \hspace{2mm} s\in {\mathbb R}^{+},
\end{eqnarray}
for some $C_0 \in {\mathbb R}$. We shall use the shorthand notation $\phi_n=\phi\left( \frac{\mid n\mid}{M}\right)$.
We now multiply  (\ref{latg1}), with the function $\phi_n\theta_n\bar{u}_n$, $n\in {\mathbb Z}$, and we sum over all sites and keep the imaginary part. The resulting equation is
\begin{eqnarray}
\label{com1}
\frac{1}{2}\frac{d}{dt}\sum_{n\in\mathbb{Z}}\phi_n\theta_n|u_n|^2-\hat{\alpha}\mathrm{L}_1(u_n)-\hat{\beta} \mathrm{L}_2(u_n)+\hat{\delta}\sum_{n\in\mathbb{Z}}\theta_n|u_n|^2
+\hat{\zeta}\sum_{n\in\mathbb{Z}}\phi_n\theta_n|u_n|^{p+1}=\mathrm{Im}\sum_{n\in\mathbb{Z}}\phi_n\theta_ng_n\overline{u}_n,
\end{eqnarray}
where the terms $\mathrm{L}_1,\mathrm{L}_2$ are found to be
\begin{eqnarray*}
\label{com2}
\mathrm{L}_1(u_n)
&=&\sum_{n\in\mathbb{Z}}(\phi_{n+1}\theta_{n+1}-\phi_n\theta_n)(u_{1,n}u_{2,n+1}-u_{2,n}u_{1,n+1}),\;\;\;\;\;\;\;\;\\
\label{com3}
\mathrm{L}_2(u_n)
&=&2\sum_{n\in\mathbb{Z}}\phi_n\theta_n(u_{1,n}^2+u_{2,n}^2)
-\sum_{n\in\mathbb{Z}}(\phi_{n+1}\theta_{n+1}+\phi_n\theta_n)(u_{1,n}u_{1,n+1}+u_{2,n}u_{2,n+1}).
\end{eqnarray*}
We rewrite $\phi_{n+1}\theta_{n+1}-\phi_n\theta_n=(\phi_{n+1}-\phi_n)\theta_{n+1}+\phi_n(\theta_{n+1}-\theta_n)$. We then get the estimate
\begin{eqnarray}
\label{com4}
|\mathrm{L}_1(u_n)|&\leq&\sum_{n\in\mathbb{Z}}|\phi_{n+1}-\phi_n|\theta_{n+1}|u_{1,n}u_{2,n+1}-u_{2,n}u_{1,n+1}|
+\sum_{n\in\mathbb{Z}}|\theta_{n+1}-\theta_n|\phi_n|u_{1,n}u_{2,n+1}-u_{2,n}u_{1,n+1}|\nonumber\\
&\leq&
\sum_{n\in\mathbb{Z}}|\phi_{n+1}-\phi_n|\theta_{n+1}|u_{1,n}u_{2,n+1}-u_{2,n}u_{1,n+1}|
+D\sum_{n\in\mathbb{Z}}\phi_n\theta_n|u_{1,n}u_{2,n+1}|+D\sum_{n\in\mathbb{Z}}\phi_n\theta_n|u_{2,n}u_{1,n+1}|\nonumber\\
&\leq&
\frac{C\rho_1^2}{M}+D(\left(\sum_{n\in\mathbb{Z}}\phi_n\theta_n|u_{1,n}|^2\right)^{1/2}
(\left(\sum_{n\in\mathbb{Z}}\phi_n\theta_n|u_{2,n+1}|^2\right)^{1/2}\nonumber\\
&&+D(\left(\sum_{n\in\mathbb{Z}}\phi_n\theta_n|u_{2,n}|^2\right)^{1/2}
(\left(\sum_{n\in\mathbb{Z}}\phi_n\theta_n|u_{1,n+1}|^2\right)^{1/2}\nonumber\\
&\leq&
\frac{C\rho_1^2}{M}+D\underline{d}^{-1/2}(\left(\sum_{n\in\mathbb{Z}}\phi_n\theta_n|u_{1,n}|^2\right)^{1/2}
(\left(\sum_{n\in\mathbb{Z}}\phi_n\theta_{n+1}|u_{2,n+1}|^2\right)^{1/2}\nonumber\\
&&+D\underline{d}^{-1/2}(\left(\sum_{n\in\mathbb{Z}}\phi_n\theta_n|u_{2,n}|^2\right)^{1/2}
(\left(\sum_{n\in\mathbb{Z}}\phi_n\theta_{n+1}|u_{1,n+1}|^2\right)^{1/2}\nonumber\\
&\leq&
\frac{C\rho_1^2}{M}+\frac{1}{2}D\underline{d}^{-1/2}\left\{\sum_{n\in\mathbb{Z}}\phi_n\theta_n|u_n|^2+
\sum_{n\in\mathbb{Z}}\phi_n\theta_{n+1}|u_{n+1}|^2\right\}\nonumber\\
&=&
\frac{C\rho_1^2}{M}+\frac{1}{2}D\underline{d}^{-1/2}\left\{\sum_{n\in\mathbb{Z}}\phi_n\theta_n|u_n|^2+
\sum_{n\in\mathbb{Z}}\phi_{n+1}\theta_{n+1}|u_{n+1}|^2-\sum_{n\in\mathbb{Z}}(\phi_{n+1}-\phi_n)\theta_{n+1}|u_{n+1}|^2\right\}\nonumber\\
&=&
\frac{C\rho_1^2}{M}+\frac{1}{2}D\underline{d}^{-1/2}\left\{2\sum_{n\in\mathbb{Z}}\phi_n\theta_n|u_n|^2-\sum_{n\in\mathbb{Z}}(\phi_{n+1}-\phi_n)\theta_{n+1}|u_{n+1}|^2\right\}\nonumber\\
&\leq&
\frac{C\rho_1^2}{M}+D\underline{d}^{-1/2}\sum_{n\in\mathbb{Z}}\phi_n\theta_n|u_n|^2.
\end{eqnarray}
Now for the second term of the rhs of $\mathrm{L}_2(u_n)$, we have
\begin{eqnarray}
\label{com5}
&&\left|\sum_{n\in\mathbb{Z}}(\phi_{n+1}\theta_{n+1}+\phi_n\theta_{n})(u_{1,n}u_{1,n+1}+u_{2,n}u_{2,n+1})\right|
\leq\sum_{n\in\mathbb{Z}}\phi_{n+1}\theta_{n+1}|u_{1,n}u_{1,n+1}|+\sum_{n\in\mathbb{Z}}\phi_n\theta_{n}|u_{1,n}u_{1,n+1}|\nonumber\\
&&+\sum_{n\in\mathbb{Z}}\phi_{n+1}\theta_{n+1}|u_{2,n}u_{2,n+1}|+\sum_{n\in\mathbb{Z}}\phi_n\theta_{n}|u_{2,n}u_{2,n+1}|\nonumber\\
&\leq&
\left(\sum_{n\in\mathbb{Z}}\phi_{n+1}\theta_{n+1}|u_{1,n}|^2\right)^{1/2}
\left(\sum_{n\in\mathbb{Z}}\phi_{n+1}\theta_{n+1}|u_{1,n+1}|^2\right)^{1/2}
+
\left(\sum_{n\in\mathbb{Z}}\phi_n\theta_{n}|u_{1,n}|^2\right)^{1/2}
\left(\sum_{n\in\mathbb{Z}}\phi_n\theta_{n}|u_{1,n+1}|^2\right)^{1/2}\nonumber\\
&&+
\left(\sum_{n\in\mathbb{Z}}\phi_{n+1}\theta_{n+1}|u_{2,n}|^2\right)^{1/2}
\left(\sum_{n\in\mathbb{Z}}\phi_{n+1}\theta_{n+1}|u_{2,n+1}|^2\right)^{1/2}
+
\left(\sum_{n\in\mathbb{Z}}\phi_n\theta_{n}|u_{2,n}|^2\right)^{1/2}
\left(\sum_{n\in\mathbb{Z}}\phi_n\theta_{n}|u_{2,n+1}|^2\right)^{1/2}\nonumber\\
&\leq&
\frac{1}{2}\left\{\sum_{n\in\mathbb{Z}}\phi_{n+1}\theta_{n+1}|u_n|^2
+\sum_{n\in\mathbb{Z}}\phi_{n+1}\theta_{n+1}|u_{n+1}|^2
+\sum_{n\in\mathbb{Z}}\phi_{n}\theta_{n}|u_{n}|^2
+\sum_{n\in\mathbb{Z}}\phi_{n}\theta_{n}|u_{n+1}|^2\right\}\nonumber\\
&\leq&
\frac{1}{2}\left\{\overline{d}\sum_{n\in\mathbb{Z}}\phi_{n+1}\theta_{n}|u_n|^2
+\sum_{n\in\mathbb{Z}}\phi_{n+1}\theta_{n+1}|u_{n+1}|^2
+\sum_{n\in\mathbb{Z}}\phi_{n}\theta_{n}|u_{n}|^2
+\underline{d}^{-1}\sum_{n\in\mathbb{Z}}\phi_{n}\theta_{n+1}|u_{n+1}|^2\right\}\nonumber\\
&=&
\sum_{n\in\mathbb{Z}}\phi_{n}\theta_{n}|u_{n}|^2+\frac{\overline{d}}{2}\sum_{n\in\mathbb{Z}}\phi_{n}\theta_{n}|u_{n}|^2
+\frac{\overline{d}}{2}\sum_{n\in\mathbb{Z}}(\phi_{n+1}-\phi_{n})\theta_{n}|u_{n}|^2\nonumber\\
&&+
\frac{\underline{d}^{-1}}{2}\sum_{n\in\mathbb{Z}}\phi_{n+1}\theta_{n+1}|u_{n+1}|^2
-\frac{\underline{d}^{-1}}{2}\sum_{n\in\mathbb{Z}}(\phi_{n+1}-\phi_{n})\theta_{n+1}|u_{n+1}|^2\nonumber\\
&\leq& \left(1+\frac{\overline{d}}{2}+\frac{\underline{d}^{-1}}{2}\right)\sum_{n\in\mathbb{Z}}\phi_{n}\theta_{n}|u_{n}|^2
+\frac{C\rho_1^2}{M}.
\end{eqnarray}
Inserting (\ref{com2}) and (\ref{com5}) to (\ref{com1}), we obtain the differential inequality
\begin{eqnarray}
\frac{1}{2}\frac{d}{dt}\sum_{n\in\mathbb{Z}} \phi_n\theta_n \mid
u_n\mid^2  +\sigma_0\sum_{n\in\mathbb{Z}}\phi_n\theta_n \mid
u_n\mid^2 \le \frac{C}{M}\rho_1^2 +\frac{1}{2\epsilon}\sum_{\mid
n \mid
>M}\theta_n\mid g_n\mid^2 .\nonumber
\end{eqnarray}
Using Gronwall's inequality, we obtain the estimate
\begin{eqnarray}
\sum_{n\in\mathbb{Z}}\phi_n\theta_n \mid u_n \mid^2 \le
e^{-2\sigma_0(t-t_0)}\sum_{n\in\mathbb{Z}}\phi_n\theta_n \mid
u_n(t_0)\mid^2 +\frac{1}{2\sigma_0}\left( \frac{2C}{M}\rho_1^2
+\frac{1}{\epsilon}\sum_{\mid n \mid >M}\theta_n\mid  g_n\mid^2  \right),
\nonumber
\end{eqnarray}
for $t>t_0$, where $t_0$ is the time of entry of initial data bounded in $\ell^2_{\theta}$, into the absorbing ball of radius $\rho_1$ in $\ell^{2}_{\theta}$.
Since $g\in \ell^{2}_{\theta}$, then  for all $\eta >0$, there exists $K(\eta)$ such that
\begin{eqnarray}
\frac{2C}{M}\rho_1^2 +\frac{1}{\epsilon}\sum_{\mid n \mid >M}\theta_n \mid g_n\mid^2  \le \eta, \hspace{3mm} \forall M > K(\eta).
\nonumber
\end{eqnarray}
Therefore, for all $\eta$ and for $t>t_0$ and $M>K(\eta)$, we obtain
that
\begin{eqnarray}
\sum_{n\in\mathbb{Z}}\phi_n\theta_n \mid u_n \mid^2 \le
e^{-2\sigma_0(t-t_0)}\rho_1^2 +\frac{1}{2\sigma_0}\eta. \nonumber
\end{eqnarray}
Choosing $t$ large enough, we may then obtain
\begin{eqnarray*}
\sum_{\mid n\mid >2M}\theta_n\mid u_n\mid^2 \le
\sum_{n\in\mathbb{Z}}\phi_n\theta_n\mid u_n\mid^2 \le
\frac{\eta}{\sigma_0} .\nonumber
\end{eqnarray*}
This estimate holds as long as $t\ge T(\eta)$, where
$$T(\eta)=t_0+\frac{1}{2\sigma_0}\ln\left(\frac{2\sigma_0\rho_1^2}{\eta}\right),$$
and $M> K(\eta)$. This concludes the proof of the Lemma.
$\diamond$

With Lemma \ref{DCGLtail} at hand, we are able to  prove that the semigroup $S(t):\ell^2_{\theta}\rightarrow\ell^2_{\theta}$ is asymptotically compact. The proof follows closely that of \cite{Bates}, adapted in the case of $\ell^2_{\theta}$, and is presented for the completeness of the presentation.
\begin{rigor2}
\label{asymcomp}
The semigroup $S(t)$ is asymptotically compact in $\ell^{2}_{\theta}$, that is, if the sequence $u_n$ is bounded in $\ell^{2}_{\theta}$ and $t_n \rightarrow \infty$, then $S(t_n)u_n$ is precompact in $\ell^2_{\theta}$.
\end{rigor2}
{\bf Proof:} It follows from Lemma \ref{ballweighted}, that if $u_n \in
\ell^2_{\theta}$, such that $\mid\mid u_n\mid\mid_{\ell^2_{\theta}}\le r$, $r>0$,
there exists $T(r)>0$ and an integer $N_1(r)$, such
that  $t_n \ge T(r)$ for $n \ge N_1(r)$, and
\begin{eqnarray}
\label{prese}
S(t_n)u_n \subset {\cal B}_0, \hspace{3mm} \forall n \ge N_1(r).
\end{eqnarray}
From (\ref{prese}), there exists  $u_0 \in \mathcal{B}_0$ and a subsequence of $S(t_n)u_n$ (not relabelled), such that $S(t_n)u_n \rightharpoonup u_0$ in $\mathcal{B}_0$. 
Lemma
\ref{DCGLtail}, implies the  existence of some $K_1(\eta)$
and $T_1(\eta)$ such that
\begin{eqnarray}
\sum_{\mid i\mid \ge K_1(\eta)}\theta_{i} \mid (S(t)S(T_r)u_n)_i \mid^{2} \le \frac{\eta^2}{8}, \hspace{3mm} t \ge T_1(\eta),
\nonumber
\end{eqnarray}
where by $(S(t)u_n )_i$, we denote the $i$-th coordinate, of the infinite sequence $S(t)u_n \in \ell^{2}_{\theta}$.
Since $t_n \rightarrow \infty$, we may find $N_2(r,\eta)\in {\mathbb N}$, such that $t_n \ge T_r +T(\eta)$, if $n \ge N_2(r,\eta)$. Hence,
\begin{eqnarray}
\sum_{\mid i\mid \ge K_1(\eta)}\theta_i \mid (S(t_n)u_n)_i \mid^{2} = \sum_{\mid i\mid \ge K_1(\eta)}\theta_i\mid (S(t_n-T_r)S(T_r)u_n)_i \mid^{2} \le \frac{\eta^2}{8}.
\nonumber
\end{eqnarray}
On the other hand, 
$$ \sum_{\mid i\mid \ge K_2(\eta)}\theta_i \mid (u_0)_i \mid^2 \le \frac{\eta^2}{8}.$$
Choosing now $K(\eta)=\mathrm{max}(K_1(\eta),K_2(\eta))$, we get for all $\eta>0$, that
\begin{eqnarray}
\mid\mid S(t_n)u_n -u_0 \mid\mid_{\ell^{2}_{\theta}}^2 = \sum_{\mid i\mid \le K(\eta)}\theta_i\mid (S(t_n)u_n -u_0)_i\mid^2 +
\sum_{\mid i\mid > K(\eta)}\theta_i\mid (S(t_n)u_n -u_0)_i\mid^2
\nonumber \\
\le \frac{\eta^2}{2}+ 2 \sum_{\mid i\mid > K(\eta)}\theta_i(\mid (S(t_n)u_n)_i \mid^2 + \mid(u_0)_i\mid^2) \le \eta^{2}
\nonumber
\end{eqnarray}
Note that the first estimate, comes from the strong convergence in the finite dimensional space ${\mathbb C}^{2K(\eta)+1}$. \ $\diamond$

The main result of this section, which is a consequence of Proposition \ref{asymcomp} and 
\cite[Theorem 1.1.1]{RTem88}, can be stated as follows:
\begin{rigor1}
\label{dynamics1} 
The semigroup $S(t)$ associated to (\ref{latg1})-(\ref{latg2})
possesses a global attractor
$\mathcal{A}=\omega(\mathcal{B}_0)\subset\mathcal{B}_0\subset\ell^2_{\theta}$
which is compact, connected and maximal among the functional
invariant sets in $\ell^2_{\theta}$.
\end{rigor1}
\subsection{The finite dimensional approximation of the global attractor in $\ell^2_{\theta}$}
This section, is devoted to the finite approximation of the global attractor, of exponentially localized solutions of (\ref{latg1})-(\ref{latg2}). 
Since from Theorems \ref{lthetaex} and \ref{dynamics1}, the solution of (\ref{latg1})-(\ref{latg2}) is in $\mathrm{C}^1(\mathbb{R}^+,\ell^2_{\theta})$, an immediate consequence of the inclusion relation (\ref{we1}), is that
\begin{eqnarray*}
\lim_{n\rightarrow\infty}u_n(t)=0,\;\;t\geq 0.
\end{eqnarray*}
Thus it is natural to  seek for approximations of the global attractor, generated by the the following system of ordinary differential equations, supplemented with Dirichlet boundary conditions, 
\begin{eqnarray}
\label{latgfin}
i\dot{v}_n+(\hat{\alpha}+i\hat{\beta})(v_{n-1}&-&2v_n+v_{n+1})+(\hat{\gamma}+i\hat{\delta})v_n+(\hat{\eta}+i\hat{\zeta})F(v_n)=g_n,\;\;|n|\leq N,\\
\label{latgD1}
v_{-(N+1)}(\cdot)&=&v_{(N+1)}(\cdot)=0,\\
\label{latgD2}
v_n(0)&=&u_{n,0},\;\;|n|\leq N.
\end{eqnarray}
System (\ref{latgfin}) can be written as an evolution equation in $\mathbb{C}^{2N+1}$,  this time endowed with the inner
product and induced norm
\begin{eqnarray}
\label{Dn}
(u,v)_{2_{\theta}}:=\mathrm{Re}\sum_{n=-N}^{n=N}\theta_n u_n\overline{v_n},\;\;||u||_{2_{\theta}}:=\sum_{n=-N}^{n=N}\theta_n|u_n|^2,\;\;u,\,v\in \mathbb{C}^{2N+1},
\end{eqnarray}
Since all the norms in the finite dimensional space $\mathbb{C}^{2N+1}$ are equivalent, a result similar to Proposition \ref{fattr}, can be produced, establishing the existence of global attractor in $X_{\theta}:=(\mathbb{C}^{2m+1},\;||\cdot||_{2_{\theta}})$, with entry time independent of the initial data. However, since for the finite dimensional approximation, we are interested in a-priori bounds in $\mathbb{C}^{2N+1}$ endowed with the $||\cdot||_{2_{\theta}}$-norm, which should be independent of $N$, it is crucial to follow the procedure described in Lemma \ref{ballweighted}, and pose the same conditions on the parameters. Thus working exactly as in Lemma \ref{ballweighted}, we may prove the following
\begin{rigor2}
\label{fattrb}
Let $v_0:=(v_{n,0})_{|n|\leq N}\in X_{\theta}$. For $1<p <\infty$, there exists a unique solution of
(\ref{latgfin})-(\ref{latgD2})
such that $v\in\mathrm{C}^1([0,\infty),X_{\theta})$. Assume further that condition (\ref{parcrit0}) holds.
Then the dynamical system defined by  (\ref{latgfin})-(\ref{latgD2}),
\begin{eqnarray}
\label{fdynamicalb} S_N(t):v_0\in X_{\theta}\rightarrow
v(t)\in X_{\theta},
\end{eqnarray}
possesses a bounded absorbing set $\mathcal{O}_0$ in
$X_{\theta}$ and a global attractor
$\mathcal{A}_{N}=\omega(\mathcal{O}_0)\subset\mathcal{O}_0\subset X_{\theta}$ :
For every bounded set $\mathcal{O}$ of $X_{\theta}$, there exists
$t_1(\mathcal{O},\mathcal{O}_0)$ such that for all  $t\geq
t_1(\mathcal{O},\mathcal{O}_0)$, it holds that
$S_N(t)\mathcal{O}\subset\mathcal{O}_0$, and for every $t\geq 0$
$S_N(t)\mathcal{A}_N=\mathcal{A}_N$.
\end{rigor2}
 Following \cite{Bates,SZ2} (see also
\cite{bab90,2Nikoi-a} for a similar idea applied to pdes
considered in all of $\mathbb{R}^N$), we observe that the
$X_{\theta}$-solution of (\ref{latgfin})-(\ref{latgD2}) can be
extended naturally in the infinite dimensional space $\ell^2_{\theta}$ , as
\begin{equation}
\label{ZKGL}
(u_N(t))_{N\in\mathbb{Z}}=\left\{
\begin{array}{cc}
v(t)=(v_n(t))_{|n|\leq N},\;\; & |n|\leq N,\\
0,\;\; & |n|> N. \nonumber
\end{array}
\right.
\end{equation}
Let us note that in the light of (\ref{ZKGL}) the finite dimensional space $X_{\theta}$, could be viewed as a finite dimensional subspace of $\ell^2_{\theta}$, with elements $u\in\ell^2_{\theta}$ satisfying the Dirichlet boundary conditions (\ref{latgD1}).

The  global attractor
$\mathcal{A}$ of the semigroup $S(t)$ associated with
(\ref{latg1})-(\ref{latg2}), will be approximated by the global
attractor $\mathcal{A}_N$ of $S_N(t)$ associated to
(\ref{latgfin})-(\ref{latgD1}), as $N\rightarrow\infty$. Recall that
the semidistance of two nonempty compact subsets of a metric space
$X$, endowed with the metric $d_X(\cdot ,\cdot)$, is defined as
\begin{eqnarray*}
d(\mathcal{B}_1,\mathcal{B}_2)=\sup_{x\in\mathcal{B}_1}\inf_{y\in\mathcal{B}_2}d_X(x,y).
\end{eqnarray*}
\begin{rigor1}
\label{dynamics2}
The global attractor $\mathcal{A}_N$ converges to $\mathcal{A}$ in the sense of the semidistance related to $\ell^2_{\theta}$: we have that $\lim_{N\rightarrow\infty}d(\mathcal{A}_N,\mathcal{A})=0$.
\end{rigor1}
{\bf Proof:}\ \ We denote by $\mathcal{U}$ an open-neighborhood of
the absorbing ball $\mathcal{B}_0$ of $S(t)$. Obviously
$\mathcal{A}$ attracts $\mathcal{U}$. For arbitrary
$N\in\mathbb{N}$, we consider the  semigroup $S_N(t)$ defined by
Proposition \ref{fattrb} and its global attractor $\mathcal{A}_N$.
Exactly as in Lemma \ref{ballweighted}, it can be shown that
$\mathcal{B}_0\cap X_{\theta}$ is also an absorbing set
for $S_N(t)$. Therefore
$$\mathcal{A}_N\subset\mathcal{B}_0\cap X_{\theta}\subset\mathcal{U}\cap X_{\theta},$$
and $\mathcal{A}_N$ attracts $\mathcal{U}\cap X_{\theta}$.
In the light of Proposition \ref{fattrb} and \cite[Theorem
I1.2, pg. 28]{RTem88}, it remains to verify that for every compact interval
$\mathrm{I}$ of $\mathbb{R}^+$,
\begin{eqnarray}
\label{cond2} \delta_N(\mathrm{I}):=\sup_{v_0\in\mathcal{U}\cap X_{\theta}}
\sup_{t\in\mathrm{I}}d(S_N(t)\psi_0,S(t)\psi_0)\rightarrow
0,\;\;\mbox{as}\;\;N\rightarrow\infty .
\end{eqnarray}
We consider the corresponding solution $v(t)=S_N(t)v_0$,
$v(0)=v_0$, in $X_{\theta}$ through
(\ref{latgfin})-(\ref{latgD1}). Then by Proposition \ref{fattr},  it
follows that $v(t)\in\mathcal{A}_N$ for any $t\in\mathbb{R}^+$.
Therefore, if $\rho>0$ is the $N$-independent radius of the
absorbing ball $\mathcal{O}_N$ in  $X_{\theta}$, then
for every $t\in\mathbb{R}^+$, $||v(t)||^2_{2_{\theta}}\leq\rho^2$.  
Using (\ref{ZKGL}), we may construct the extension of $v(t)$ in $\ell^2_{\theta}$. The extension $u_N(t)$
satisfies the estimates
\begin{eqnarray}
\label{aprrox1}
||u_N(t)||_{\ell^2_{\theta}}^2\leq\rho^2,\;\;||\dot{u}_N(t)||_{\ell^2_{\theta}}^2\leq
C(\rho, ||g||_{\ell^2_{\theta}}),
\end{eqnarray}
the latter derived by (\ref{latgfin}).
According to \cite[Theorem 10.1 pg. 331-332]{RTem88} or \cite[Lemma 4, pg. 60]{SZ2}, for the  justification of (\ref{cond2}) it suffices to show that $u_N(t)$ converges to  a solution $u(t)$ of (\ref{latg1})-(\ref{latg2}) in an arbitrary compact interval of $\mathbb{R}^+$, and $u_0=u(0)$ in a bounded set of $\ell^2_{\theta}$. Let $\mathrm{I}$ be an arbitrary compact interval of $\mathbb{R}^+$. 
From estimates (\ref{aprrox1}), we may extract a subsequence $u_j$ of $u_N$, such that
\begin{eqnarray}
\label{dweakc} u_{j}(t)\rightharpoonup
u(t),\;\;\mbox{in}\;\;\ell^2_{\theta},\;\;\mbox{as}\;\;
j\rightarrow\infty,\;\;\mbox{for every}\;\;t\in D,
\end{eqnarray}
where $D$ denotes a countable dense subset of $\mathrm{I}$.

For any $t\in \mathrm{I}$ we consider the  sequence
$$\chi_N(t):=(u_N(t),z)_{\ell^2_{\theta}},\;\;z\in\ell^2_{\theta},$$
which by (\ref{aprrox1}), is differentiable as a function of $t$, and 
$\chi_N'(t)=(\dot{u}_N(t),z)_{\ell^2_{\theta}}$. Moreover,
there exists $\xi\in \mathrm{I}$ such that, for fixed $t,s\in
\mathrm{I}$
\begin{eqnarray}
\label{equco1}
|\chi_N(t)-\chi_N(s)|=|(u_N(t)-u_N(s),z)_{\ell^2_{\theta}}|&=&|(\dot{u}_N(\xi),z)_{\ell^2_{\theta}}|\,|t-s|\nonumber\\
&\leq&\sup_{\xi\in
I}||\dot{u}_N(\xi)||_{\ell^2_{\theta}}||z||_{\ell^2_{\theta}}|t-s|\leq C|t-s|.
\end{eqnarray}
i.e the sequence $\chi_N$ is equicontinuous. On the other hand it follows from (\ref{equco1}), that there exists an $N$-independent constant $C_1$ such that
\begin{eqnarray}
\label{Ascoli} ||u_N(t)-u_N(s)||_{\ell^2}\leq C_1|t-s|,
\end{eqnarray}
Hence by Ascoli's Theorem, it follows that the convergence (\ref{dweakc}), holds uniformly on $\mathrm{I}$ as $N\rightarrow\infty$. Summarizing, we obtain  for the 
subsequence $u_j$,  the convergence relations
\begin{eqnarray}
\label{passlim}
&&u_{j}\rightarrow u \;\;\mbox{in}\;\;\mathrm{C}(\mathrm{I},\ell^2_{\theta}),\nonumber\\
&&u_{j}\stackrel{*}{\rightharpoonup} u\;\;\mbox{in}\;\;\mathrm{L}^{\infty}(\mathrm{I},\ell^2_{\theta}),\\
&&\dot{u}_{j}\stackrel{*}{\rightharpoonup}
\dot{u}\;\;\mbox{in}\;\;\mathrm{L}^{\infty}(\mathrm{I},\ell^2_{\theta}).\nonumber
\end{eqnarray}

For the passage to the limit, we shall use an eqivalent formulation of (\ref{latg1})-(\ref{latg2}). Clearly, any solution of (\ref{latg1})-(\ref{latg2}), satisfies for every $v\in\ell^2_{\theta}$ and $z (t)\in C^{\infty}_0(\mathrm{I})$, the formula 
\begin{eqnarray}
\label{weakform}
\int_{\mathrm{I}}(i\dot{u}(t),v)_{\ell^2_{\theta}}z (t)dt +
\int_{\mathrm{I}}(\mathbf{L}u(t),v)_{\ell^2_{\theta}}z (t)dt
&+&\int_{\mathrm{I}}((\hat{\gamma}+i\hat{\delta})u(t),v)_{\ell^2_{\theta}}z(t)dt\nonumber\\
&+&\int_{\mathrm{I}}((\hat{\eta}+i\hat{\zeta})F(u(t)),v)_{\ell^2_{\theta}}z(t) dt
=\int_{\mathrm{I}}(g,v)_{\ell^2_{\theta}}z(t) dt,
\end{eqnarray}
where $F(u)=|u|^{p-1}u$. Since for fixed $N\in\mathbb{Z}^+$, $u_N$ is a solution of (\ref{latgfin})-(\ref{latgD2}), we may reproduce (\ref{weakform}) for $u_{N}$, by multiplying $(\ref{latgfin})$ by $v\in\ell^2_{\theta}$, in the $\ell^2_{\theta}$-scalar product.  
By Lemma \ref{LipN}, $F:\ell^2_{\theta}\rightarrow\ell^2_{\theta}$ is
Lipschitz continuous on bounded sets of $\ell^2_{\theta}$. Therefore from (\ref{aprrox1}), it follows that there exists
a constant $c(\rho)$, such that
$||F(u_j)-F(u)||_{\ell^2_{\theta}}\leq
c(\rho)||u_{j}-u||_{\ell^2_{\theta}}$. Then from (\ref{passlim}), we
infer
\begin{eqnarray*}
\label{pasnon1}
\left|\int_{\mathrm{I}}(F(u_j(t))-F(u),v)_{\ell^2_{\theta}}z(t)dt\right|
&\leq&\int_{\mathrm{I}}||F(u_j(t))-F(u(t))||_{\ell^2_{\theta}}||v||_{\ell^2}z(t)dt\\
&\leq&c\int_{\mathrm{I}}||u_{j}(t)-u(t)||_{\ell^2_{\theta}}||v||_{\ell^2_{\theta}}z (t)\\
&\leq&c\sup_{t\in\mathrm{I}}||u_j(t)-u(t)||_{\ell^2_{\theta}}||v||_{\ell^2_{\theta}}\int_{\mathrm{I}}|z(t)|dt\rightarrow
0,\;\;\mbox{as}\;\;\mu\rightarrow\infty.
\end{eqnarray*}
On the other hand, we have from Lemma \ref{thm:LOCALLIPS2},
that $\mathbf{L} :\ell^2_{\theta}\rightarrow \ell^2_{\theta}$
is globally Lipschitz on $\ell^2_{\theta}$. Hence, we have that
\begin{eqnarray*}
\label{pasnon2}
\left|\int_{\mathrm{I}}(\mathbf{L}u_j(t)-\mathbf{L}u(t),v)_{\ell^2_{\theta}}z(t)dt\right|
&\leq&\int_{\mathrm{I}}||\mathbf{L}u_j(t)-\mathbf{L}u(t)||_{\ell^2_{\theta}}||v||_{\ell^2}z(t)dt\\
&\leq&L\int_{\mathrm{I}}||u_{j}(t)-u(t)||_{\ell^2_{\theta}}||v||_{\ell^2_{\theta}}z (t)\\
&\leq&L\sup_{t\in\mathrm{I}}||u_j(t)-u(t)||_{\ell^2_{\theta}}||v||_{\ell^2_{\theta}}\int_{\mathrm{I}}|z(t)|dt\rightarrow
0,\;\;\mbox{as}\;\;\mu\rightarrow\infty.
\end{eqnarray*}
Since $I$ is arbitrary, (\ref{weakform}) is satisfied for all $t\in\mathbb{R^+}$, i.e. $u(t)$ solves (\ref{latg1})-(\ref{latg2}). Moreover by (\ref{passlim}), we get that 
$u(t)$ is bounded in $\ell^2_{\theta}$ for all $t\in\mathbb{R^+}$. Therefore $u(t)\in\mathcal{A}$, which implies that $u_{j}(0)\rightarrow u(0)$, and $u(0)$ is at least, in a bounded set of $\ell^2_{\theta}$.   Since the convergence holds for any other subsequence having the above formulated properties, by a contradiction argument using uniqueness, we may deduce that the convergence holds for the original sequence $u_N$. Condition
(\ref{cond2}) is proved.\ \ $\diamond$.
\begin{remark} {\bf The global attractor of exponentially localized solutions for the DCGL and DNLS equations.}\ 
{\em To recover from the complex lattice differential equation (\ref{latg1}),  the DCGL equation 
\begin{eqnarray*}
\dot{u}_n-(\lambda+i\alpha)(u_{n-1}&-&2u_n+u_{n+1})+\gamma u_n+(k+i\beta)|u_n|^{p-1}u_n=f_n,\;\;1<p<\infty,\\
u_n(0)&=&u_{n,0},
\end{eqnarray*}
(with gauge nonlinear interaction $(\mathcal{G})$-note that this time the nonlinearity is on the lhs of the equation ), we have to set
$\hat{\alpha}=\alpha,\hat{\beta}=-\lambda,\hat{\gamma}=0,\hat{\delta}=-\gamma,\hat{\eta}=-\beta, \hat{\zeta}=k$ and
$g_n=if_n$.
For these values of the parameters, condition  (\ref{parcrit0}) becomes
$$-\gamma-\frac{\epsilon}{2}+2\lambda-|\alpha|D\underline{d}^{1/2}-\lambda\left(1+\frac{\overline{d}}{2}+\frac{\underline{d}^{-1}}{2}\right),\;\;k>0.$$
In particular, in the case of the exponential weight $\theta_n=\exp(\mu|n|)$ (exponential localization), we find 
$\overline{d}=e^\mu$, $\underline{d}=e^{-\mu}$, $D=e^\mu-1$. We set $\epsilon=2\lambda >0$
and condition (\ref{parcrit0}) reads as
$$-\gamma > \lambda e^\mu+2|\alpha|\sinh(\mu/2),\;\;\gamma <0.$$
Note that in the absence of the external excitation ($f_n=0$), the dynamical system exhibits trivial dynamics, in the sense that  $\limsup_{t\rightarrow\infty}||u(t)||_{\ell^2_{\theta}}^2=0$.

For the weakly damped and driven DNLS, ($\lambda=k=0$)
\begin{eqnarray*} 
\dot{u}_n-i\alpha(u_{n-1}&-&2u_n+u_{n+1}) +i\gamma u_n+\beta|u_n|^{p-1}u_n =f_n,\;\;n\in\mathbb{Z},\;\;1< p <\infty
\end{eqnarray*}
the condition on the dissipation parameter is
$$-\gamma > 2|\alpha|\sinh(\mu/2),\;\;\gamma <0.$$
Although this condition appears from the consideration of the problem in $\ell^2_{\theta}$, it seems to be in conformity with the analysis of section 3.3, since there is not any restriction for the sign of the parameter $\beta$. Let us note that even for the weakly damped and undriven NLS partial differential equation there exist initial data for which solutions may blow-up in finite time see \cite{Tsu1}. This is never the case for the damped and undriven  ($f_n=0$) DNLS, for which solutions always exhibit energy decay.} 
\end{remark}

$^\dag$Department of Mathematics, \\
University of the Aegean,\\
Karlovassi, 83200 Samos, GREECE\\
E-mail address: \ \ karan@aegean.gr\\
\\
$^*$Department of Telecommunications Science and Technology\\ 
University of the Peloponesse,\\ 
Tripolis 22100, GREECE\\
E-mail address:\ \ enistaz@aegean.gr\\
\\
$^\ddag$Department of Statistics and Actuarial Science, \\
University of the Aegean,\\
Karlovassi, 83200 Samos, GREECE\\
E-mail address: \ \ ayannaco@aegean.gr
\end{document}